\documentclass[11pt]{amsart}

\usepackage[all]{xy}
\usepackage{amscd}
\usepackage{mathrsfs}
\usepackage{amssymb}
\usepackage{amsthm}
\usepackage{amsmath}

\newtheorem{theorem}{Theorem}[section]
\newtheorem{lemma}[theorem]{Lemma}
\newtheorem{proposition}[theorem]{Proposition}
\newtheorem{corollary}[theorem]{Corollary}
\newtheorem{remark}[theorem]{Remark}

\theoremstyle{definition}
\newtheorem{definition}[theorem]{Definition}
\newtheorem{example}[theorem]{Example}

\theoremstyle{remark}

\numberwithin{equation}{section}

\newcommand{\bde}{\begin{definition}}
\newcommand{\ede}{\end{definition}}
\newcommand{\bpr}{\begin{proposition}}
\newcommand{\epr}{\end{proposition}}
\newcommand{\bth}{\begin{theorem}}
\newcommand{\ethm}{\end{theorem}}
\newcommand{\bexa}{\begin{example}}
\newcommand{\eexa}{\end{example}}
\newcommand{\bcor}{\begin{corollary}}
\newcommand{\ecor}{\end{corollary}}
\newcommand{\blem}{\begin{lemma}}
\newcommand{\elem}{\end{lemma}}
\newcommand{\brem}{\begin{remark}}
\newcommand{\erem}{\end{remark}}
\newcommand{\bprf}{\begin{proof}}
\newcommand{\eprf}{\end{proof}}
\def \benu{\begin{enumerate}\renewcommand{\labelenumi}{(\roman{enumi})}\renewcommand{\itemsep}{0pt}}
\newcommand{\eenu}{\end{enumerate}}
\newcommand{\A}{\mathbb{A}}
\newcommand{\fkm}{{\mathfrak m}}

\newcommand{\fkp}{{\mathfrak p}}

\newcommand{\mC}{{\mathscr C}}

\newcommand{\mF}{{\mathscr F}}
\newcommand{\mG}{{\mathscr G}}

\newcommand{\mO}{\mathcal O}
\newcommand{\mI}{{\mathscr I}}
\newcommand{\mN}{{\mathscr N}}

\newcommand{\Z}{\mathbb{Z}}

\newcommand{\Q}{\mathbb{Q}}
\newcommand{\F}{\mathbb{F}}
\newcommand{\G}{\mathbb{G}}
\newcommand{\bH}{\mathbb{H}}
\newcommand{\Ps}{\mathbb{P}}

\newcommand{\im}{{\rm Im\,}}

\def\Ker{{\rm Ker\,}}

\newcommand{\gal}{{\rm Gal}}

\newcommand{\ind}{{\rm Ind}}

\newcommand{\id}{{\rm id}}

\DeclareMathOperator{\Spec}{\mathrm{Spec}}
\DeclareMathOperator{\Proj}{\mathrm{Proj}}

\DeclareMathOperator{\Frac}{\mathrm{Frac}}

\newcommand{\Hom}{{\rm Hom}}

\newcommand{\End}{{\rm End}}

\newcommand{\Sh}{{\rm Sh}}

\newcommand{\bsc}{\begin{tobira}}
\newcommand{\esc}{\end{tobira}}
\newcommand{\abc}{\xymatrix{0 \ar@{>}[r] & A \ar@{>}[r]^{f} & B \ar@{>}[r]^{g} & C \ar@{>}[r] & 0}}
\newcommand{\xyz}{\xymatrix{0 \ar@{>}[r] & X \ar@{>}[r]^{f} & Y \ar@{>}[r]^{g} & Z \ar@{>}[r] & 0}}
\def\phi{\varphi}
\def\hat{\widehat}
\def\lra{\longrightarrow}
\def\Lra{\Longrightarrow}

\def\ra{\rightarrow}

\def\ol{\overline}

\def\ul{\underline}
\def\MU{\mbox{\boldmath{$\mu$}}}

\DeclareMathOperator{\chara}{\mathrm{char}}
\def\hat{\widehat}
\def\tilde{\widetilde}
\def\ur{{\rm ur}}
\def\univ{{\rm univ}}
\def\cusp{{\rm cusp}}
\def\st{{\rm st}}
\def\isom{\stackrel{\cong}{\longrightarrow}}

\begin{document}

\title[Non-Abelian Lubin-Tate Theory]{On non-abelian Lubin-Tate theory via vanishing cycles}

\author[Teruyoshi Yoshida]{Teruyoshi Yoshida}
\address{University of Cambridge, Department of Pure Mathematics and
Mathematical Statistics, Centre for Mathematical Sciences, Wilberforce
Road, Cambridge, CB3 0WB, UK}
\email{T.Yoshida@dpmms.cam.ac.uk}
\date{\today}

\begin{abstract}
We give a purely local proof, in the depth 0 case, of the result by Harris-Taylor which asserts that the local Langlands correspondence for $GL_n$ is realized in the vanishing cycle cohomology of the deformation spaces of one-dimensional formal modules of height $n$. Our proof is given by establishing the direct geometric link with the Deligne-Lusztig theory for $GL_n(\F_q)$. 
\end{abstract}

\maketitle

\tableofcontents

\section{Introduction}

Let $p$ be a prime, and $K$ be a finite extension of the $p$-adic field $\Q_p$, with the ring of integers $\mO$ and the residue field $k$ of cardinality $q$. The proof of the local Langlands correspondence for $GL_n(K)$, by Harris-Taylor \cite{HT}, was achieved by showing that the desired correspondence is realized in the $\ell$-adic vanishing cycle cohomology groups of the deformation spaces of formal $\mO$-modules of height $n$ with Drinfeld level structures (known as {\em non-abelian Lubin-Tate theory} or the conjecture of Deligne-Carayol \cite{Ca2}). As these deformation spaces occur as complete local rings of certain unitary Shimura varieties at the ``supersingular" points, they made an essential use of the fact that global Langlands correspondences are realized in the $\ell$-adic etale cohomology groups of these Shimura varieties over CM fields. In this article, we give a purely local approach to this non-abelian Lubin-Tate theory, in the special case of depth 0 or level $\fkp$, by computing the local equation of the deformation space and constructing its suitable resolution to calculate the vanishing cycle cohomology directly. We show that, in this case, the non-abelian Lubin-Tate theory for supercuspidal representations of $GL_n(K)$ is essentially equivalent to the Deligne-Lusztig theory for $GL_n$ of the residue field $k$, which realizes the cuspidal representations of $GL_n(k)$ in the $\ell$-adic cohomology groups of certain varieties over an algebraic closure $\ol{k}$ of $k$.

To state our theorems precisely, let $K,\mO,k$ as above and fix $n\geq 1$. Let $K^\ur$ be the maximal unramified extension of $K$, and let $W$ be the completion of the ring of integers $\mO^\ur$ of $K^\ur$ (sometimes denoted $W_\mO(\ol{k})$ in the literature). Let $\eta,\ol{\eta}$ be the spectra of $\Frac W$ and its fixed algebraic closure, respectively.

Firstly, let $X$ be the spectrum of the deformation ring of formal $\mO$-module of height $n$ with level $\fkp$ structure (\cite{Dr}), which is a scheme of relative dimension $n-1$ over $W$. We are interested in the $\ell$-adic etale cohomology groups $H^i(X_{\ol{\eta}},\ol{\Q}_\ell)$ ($\ell\neq \chara k$) of the geometric generic fiber $X_{\ol{\eta}}:=X\times_{\Spec W}\ol{\eta}$, which are finite dimensional $GL_n(k)\times I_K$-modules, where $I_K$ is the inertia group of $K$. Secondly, let $DL$ be the Deligne-Lusztig variety for $GL_n(k)$, associated to the element of the Weyl group of $GL_n$ that corresponds to the cyclic permutation $(1,\ldots,n)$ in the symmetric group of $n$ letters, or equivalently to a non-split torus $T$ with $T(k)\cong k_n^\times$ where $k_n$ is the extension of $k$ of degree $n$ (\cite{DL}). This $DL$ is a smooth affine variety over $\ol{k}$ with actions of $GL_n(k)$ and $T(k)\cong k_n^\times$, hence we can regard $H^i_c(DL,\ol{\Q}_\ell)$ as a $GL_n(k)\times I_K$-module by the canonical surjection $I_K\ra k_n^\times$. 

We denote the alternating sums of these cohomology groups as follows:
\[H^*(X_{\ol{\eta}}):=\sum_i(-1)^i[H^i(X_{\ol{\eta}},\ol{\Q}_\ell)],\ \ \ H^*_c(DL):=\sum_i(-1)^i[H^i_c(DL,\ol{\Q}_\ell)],\]
which are regarded as elements of the Grothendieck group of $GL_n(k)\times I_K$-modules. Then our main theorem on the vanishing cycle cohomology groups of $X$ can be stated as follows (Theorem \ref{cohomology}):

\bth
\label{main1}
\benu
\item We have the equality $H^*(X_{\ol{\eta}})=H^*_c(DL)$.
\item Among the $H^i(X_{\ol{\eta}},\ol{\Q}_\ell)$, cuspidal representations $\pi$ of $GL_n(k)$ and generic inertia characters $\chi$ of $I_K$ (here {\rm generic} means $\chi$ does not factor through any $k_m^\times$ with $m\mid n,\ m<n$ via the norm map $k_n^\times \ra k_m^\times$) occur only in $H^{n-1}(X_{\ol{\eta}},\ol{\Q}_\ell)$, where they are coupled as $\bigoplus \pi_\chi\otimes\chi$ by the Deligne-Lusztig correspondence $\chi\leftrightarrow \pi_\chi$ characterized by
\[\pi_\chi \otimes {\rm St}=\ind_{T(k)}^{GL_n(k)}\chi,\]
where ${\rm St}$ is the Steinberg representation of $GL_n(k)$.
\eenu
\ethm

The correspondence in the part (ii) can essentially be deduced from one of the main theorems of Harris-Taylor \cite{HT} (Theorem VII.1.5), which was proven via highly nontrivial global arguments, but in this article we first prove the part (i) of the above theorem by a local geometric argument, and apply the results of Deligne-Lusztig theory. The fact that the supercuspidal representations appear only in the degree $n-1$ (in the limit of cohomology groups for all levels $\fkp^m$) was remarked by Faltings \cite{Fa1}, as a refinement of the results of Harris-Taylor, and recently proved by Mieda \cite{Mi} via purely local argument.

To see that the results of Harris-Taylor imply the above theorem, one only needs to spell out the depth 0 case of the local Langlands correspondence. Here, the supercuspidal representations of $GL_n(K)$ are obtained as compact inductions of the pull back of cuspidal representations of $GL_n(k)$ to $GL_n(\mO)$. The irreducible $n$-dimensional representations of the Weil group $W_K$ are obtained by extending the generic tame inertia characters $\chi$ as in the theorem from $I_K$ to the Weil group $W_L$ of the unramified extension $L$ of degree $n$ over $K$, and then inducing them from $W_L$ to $W_K$. The local Langlands correspondence (up to twists by unramified characters) boils down to the correspondence mentioned in (ii) of the above theorem. This case gives all the {\em tamely ramified} irreducible representations of $W_K$.

To prove the above theorem by purely local arguments, we construct a suitable model of the deformation space $X$ and compute the cohomology of the geometric generic fiber $X_{\ol{\eta}}$ in terms of vanishing cycle sheaves on the special fiber. In its course we obtain important information concerning the geometry of $X$ as the following:

\bth
\label{main2}
Let $\varpi$ be a uniformizer of $\mO$.
\benu
\item (Prop.\ \ref{localeq}) The $W$-scheme $X$ is isomorphic to 
\[\Spec W[[X_1,\ldots,X_n]]/(P(X_1,\ldots,X_n)-\varpi),\]
where $P\in W[[X_1,\ldots,X_n]]$ is of the form:
\[(\text{\rm unit})\cdot \!\!\!\prod_{(a_i \bmod \fkp)_i\in k^n\setminus \{\ul{0}\}}\!\!\!\bigl([a_1](X_1)+_{\widetilde{\Sigma}}\cdots +_{\widetilde{\Sigma}}[a_n](X_n)\bigr)\]
where $[a_i]$ and $+_{\widetilde{\Sigma}}$ denote the formal $\mO$-multiplication and addition of a formal $\mO$-module over $W[[X_1,\ldots,X_n]]$ obtained by lifting the universal formal $\mO$-module over $X$.
\item  (Theorem \ref{gsstmodel}) There exists a generalized semistable model $Z_{st}$ of $X$ over $W$, i.e.\ a proper $W$-morphism $Z_{st}\ra X$ which is an isomorphism on the generic fibers and $Z_{st}$ being generalized semistable. Here {\em generalized semistable} means that its complete local rings at all the closed points are isomorphic over $W$ to 
\[W[[T_1,\ldots,T_n]]/(T_1^{e_1}\cdots T_d^{e_d}-\varpi)\ \ (d\leq n),\]
where the integers $e_i\ (1\leq i\leq d)$ are all prime to $\chara k$.
\item (Prop.\ \ref{undl}) Over the tamely ramified extension $W_n:=W(\varpi^{1/(q^n-1)})$ of $W$, there is a model of $X$ whose special fiber contains a smooth affine variety over $\ol{k}$ which is isomorphic to $DL$ as schemes with right $GL_n(k)\times I_K$-action.
\eenu
\ethm

The part (i) of this theorem gives the integral local equations of the relevant unitary Shimura varieties at supersingular points (similar equations can be given for the deformation spaces of formal $\mO$-modules with level $\fkp^m$ structures for any $m\geq 1$), and in the special case $K=\Q_p$ and $n=2$, it gives the ``integral'' version of Katz-Mazur's description of the bad reduction of modular curves $X(p^m)$ (\cite{KM} Theorem 13.8.4). The resolution constructed in the proof of the part (ii) of this theorem can be used to give a generalized semistable model of the unitary Shimura variety with level $\fkp$ structure, which suggests a more ``local" approach for computing the cohomology of these Shimura varieties. The part (iii) of the theorem is the basis of the proof of Theorem \ref{main1} above, and is obtained by normalizing the base change of a relevant part of generalized semistable model to $W_n$.

The vanishing cycle cohomology of the coverings of Lubin-Tate spaces is known to incorporate the {\em local Jacquet-Langlands correspondence} as well as local Langlands correspondence (\cite{HT}). The realization of Jacquet-Langlands correspondence was proved via local arguments by M.\ Strauch (\cite{St1}, \cite{St2}, \cite{St3}) using the period map and the trace formula on rigid analytic spaces. We hope to clarify how our work is related to other works in the field (\cite{Bo2}, \cite{Da}, \cite{Fa1}, \cite{FGL}, \cite{Mi}) in the near future.

{\bf Acknowledgements.} This work was the PhD thesis of the author, submitted to the University of Tokyo in January 2004. The May 2004 version of this article was accepted for publication in the Annales de l'Institut Fourier, but the author was unable to make the necessary revision for a long time. One of the reasons was that the author wishes to present a more canonical and moduli-theoretic approach to the results in this article in the near future. However, as the explicit calculations presented in this article may be of some interest, here we left this article close to its original form. The author heartfully thanks the editor of this volume for the acceptance of this article.

The author would like to thank his adviser Takeshi Saito for his encouragements and many improvements of the proofs. The author is also grateful to his previous adviser Kazuya Kato for his warm encouragements during the graduate program. The problem treated here grew from the numerous discussions with Tetsushi Ito. Also this work is a generalization of some computations in the modular curve case ($n=2,\ K=\Q_p$) done by B.\ Edixhoven in his unpublished manuscript \cite{Ed}, whom author would like to thank sincerely for showing him the manuscript. The author would like to thank G.\ Lusztig and Y.\ Mieda for their assistance in Deligne-Lusztig theory and etale cohomology, K.\ Ban, S.\ DeBacker, M.\ Strauch for helpful comments and discussions, Y.\ Yasufuku for the manuscript \cite{Ya} from which he learned many things about formal groups. The author is indebted to the referee of this article for his very careful reading and some essential corrections. Finally, a large part of this work was done at Harvard University, and the author would like to thank Richard Taylor for many invaluable discussions and constant encouragement during the development of this article.

{\bf Notation.} For $j\geq 1$, we denote the group of $j$-th roots of unity by $\MU_j$. Cohomology groups are all $\ell$-adic etale cohomology groups, where we fix a prime $\ell$ different from the residue characteristic $p$. For a representation $V$ of a group, we denote the corresponding element in a suitable Grothendieck group by $[V]$. For a ring $A$, we denote its group of units by $A^\times$. For a field $F$, we denote a fixed separable closure of $F$ by $\ol{F}$. For a finite field $k\cong \F_q$ with $q$ elements and $n\geq 1$, we denote by $k_n\cong \F_{q^n}$ the unique extension of $k$ of degree $n$. For a scheme $X$ and its point $x\in X$ (resp.\ a geometric point $x$ of $X$), we denote the Zariski local ring (resp.\ strict local ring) at $x$ by $\mO_{X,x}$ or $\mO_x$. Sometimes we refer to the elements of the coordinate rings as ``coordinates'', but when there is a risk of confusion we distinguish them by capital/lower cases.

\section{The deformation spaces of formal $\mO$-modules}

Let $p$ be a fixed prime, and $K$ be a finite extension of the $p$-adic field $\Q_p$. Equivalently, it is a complete discrete valuation field of characteristic zero with a finite residue field. Let $\mO$ be the ring of integers of $K$, and $\fkp$ be its maximal ideal. Let $q$ be the cardinality of the residue field, which we denote by $k:=\mO/\fkp\cong \F_q$. These are fixed throughout the paper, with the exception of Section 6.1.

We denote the maximal unramified extension of $K$ by $K^\ur$, and its completion by $\hat{K}:=\hat{K}^\ur$. We denote the ring of integers of $\hat{K}$ by $W$, which is a complete discrete valuation ring with the residue field $\ol{k}=\ol{\F}_q$. The uniformizers of $\mO$ are also uniformizers of $W$. We often fix a uniformizer $\varpi$ of $\mO$.

Here we recall from \cite{Dr} the basic facts on the deformation theory of one-dimensional formal $\mO$-modules. The author is greatly indebted to the expository article of Yasufuku \cite{Ya}, which gives a detailed account of Drinfeld's theory.

\subsection{Formal $\mO$-modules}

For a (commutative) $\mO$-algebra $A$, by a {\em formal $\mO$-module} over $A$, we mean a pair $\Sigma=(F,[\cdot])$, where $F(X,Y)\in A[[X,Y]]$ is a one-dimensional commutative formal group law and $[\cdot]:\mO\ni a\mapsto [a]\in \End(F)$ is an injective ring homomorphism such that $[a](X)\equiv aX\pmod{X^2}$. Here $\End(F):=\Hom(F,F)$, where
\[ \Hom(F,G):= \bigl\{ f(X)\in X\cdot A[[X]] \bigm| f(F(X,Y))=G(f(X),f(Y)) \bigr\} \]
for formal group laws $F,G$ over $A$. The homomorphisms between formal $\mO$-modules $\Sigma=(F,[\cdot]),\ \Sigma'=(F',[\cdot]')$ are the elements of $\Hom(F,F')$ which commute with $\mO$-multiplications. 
For any $A$-algebra $B$ and a formal $\mO$-module $\Sigma$ over $A$, a formal $\mO$-module $\Sigma\otimes_AB$ over $B$ is defined by the images of $F$ and $[a]$ under the induced homomorphisms $A[[X,Y]]\ra B[[X,Y]]$ and $A[[X]]\ra B[[X]]$. For $\Sigma=(F,[\cdot])$, we often use the notation:
\[ X+_\Sigma Y:=F(X,Y),\ \ [\cdot]_\Sigma:=[\cdot].\]

\bexa
The {\it additive group} $\G_a=(F,[\cdot])$, over arbitrary $A$, is defined by $F(X,Y)=X+Y$ and $[a](X)=aX$ for every $a\in \mO$. 
\eexa

\bpr 
{\rm (Drinfeld \cite{Dr})}\ For a formal $\mO$-module $\Sigma$ over $\ol{\F}_q$, not isomorphic to the additive group $\G_a$, there is a unique integer $n\geq 1$ (called the {\em height} of $\Sigma$) such that if $\varpi$ is a uniformizer of $\mO$ then $[\varpi]_\Sigma(X)=u\cdot X^{q^n}$ for some $u\in \ol{\F}_q[[X]]^\times$. We define the height of $\G_a$ to be $\infty$. The formal $\mO$-modules over $\ol{\F}_q$ are classified up to isomorphism by their heights.
\epr

Let $n\geq 1$, and fix a formal $\mO$-module $\Sigma_n$ of height $n$ over $\ol{\F}_q$, which is unique up to isomorphism. If $\varpi$ is a uniformizer of $\mO$, then we can (and will) choose $\Sigma_n=(F,[\cdot])$ to be {\em $\varpi$-normal}, i.e.\ it satisfies the following:
\benu
\item $[\varpi](X)=X^{q^n}$.
\item $F(X,Y)\in \F_{q^n}[[X,Y]],\ \ F(X,Y)\equiv X+Y \pmod{\deg q^n}$.
\item $[a](X)\in \F_{q^n}[[X]],\ \ [a](X)\equiv aX \pmod{X^{q^n}}$ for every $a\in \mO$. 
\eenu


\subsection{Deformations of formal $\mO$-modules}

Let $\mC$ be the category of complete noetherian local $W$-algebras $(A,\fkm_A)$ such that the structure morphisms $W\ra A$ induce isomorphisms between the residue fields $\ol{k}=W/\fkp W\isom A/\fkm_A$. We identify $\ol{k}$ and $A/\fkm_A$ by this isomorphism. The morphisms in $\mC$ are local $W$-homomorphisms. For a formal $\mO$-module $\Sigma$ over $A\in \mC$, its {\em reduction} $\bmod \fkm_A$ is defined as $\ol{\Sigma}:=\Sigma\otimes_A\ol{k}$, which is a formal $\mO$-module over $\ol{k}$.

Let us fix an integer $n\geq 1$ in the rest of the paper. For $A\in \mC$, a {\em deformation of $\Sigma_n$ to $A$} is a pair $(\Sigma,i)$, consisting of a formal $\mO$-module $\Sigma$ over $A$ and an isomorphism $i:\Sigma_n\isom \ol{\Sigma}$ as formal $\mO$-modules over $\ol{k}$. Two deformations $(\Sigma,i)$ and $(\Sigma',i')$ are {\em equivalent} if there is an isomorphism $f:\Sigma\cong \Sigma'$ such that $i'=(f \bmod \fkm_A)\circ i$. The deformation functor $\mF_0$ from $\mC$ to the category of sets is defined by sending $(A,\fkm_A)\in \mC$ to the set of equivalence classes of deformations of $\Sigma_n$ over $A$.

\bpr
{\rm (Drinfeld \cite{Dr})} The deformation functor $\mF_0$ is representable by a ring $A_0\in\mC$, which is isomorphic to the formal power series ring of $n-1$ variables $W[[T_1,\ldots,T_{n-1}]]$ over $W$. We denote the universal formal $\mO$-module over $A_0$ by $\Sigma^\univ_n$.
\epr

\subsection{Deformations with Drinfeld level structure}

For a formal $\mO$-module $\Sigma$ over $A\in \mC$, the maximal ideal $\fkm_A$ of $A$ is endowed with an $\mO$-module structure by defining the addition and the $\mO$-multiplication by 
\[x+_\Sigma y=F(x,y),\ \ [a](x)=[a]_\Sigma(x)\ \ (x,y\in \fkm_A,\ a\in \mO). \]
We denote this $\mO$-module by $\fkm_\Sigma$. A homomorphism $f:\Sigma\ra \Sigma'$ of formal $\mO$-modules over $A$ induces an $\mO$-homomorphism $f_\fkm:\fkm_\Sigma\ra \fkm_{\Sigma'}$ of the corresponding $\mO$-modules.

For a formal $\mO$-module $\Sigma=(F,[\cdot])$ over $A\in \mC$ of height $n$ and an integer $m\geq 1$, a {\it Drinfeld level $\fkp^m$ structure} on $\Sigma$ is defined to be an $\mO$-module homomorphism $\phi: (\fkp^{-m}/\mO)^n\ra \fkm_\Sigma$, which satisfies the divisibility:
\[\prod_{x\in (\fkp^{-m}/\mO)^n}\!\!\!(X-\phi(x))\ \ \bigg|\ \ [\varpi^m](X)\]
in $A[[X]]$, for some (equivalently, every) uniformizer $\varpi$ of $\mO$. 

Let us fix a uniformizer $\varpi$ of $\mO$ and choose $\Sigma_n$ to be $\varpi$-normal. If we denote the left hand side of the above by $P_\phi(X):=\prod (X-\phi(x))$, and let $[\varpi^m](X)=U(X)P_\phi(X)$ with $U(X)\in A[[X]]$, then the constant term $u_\Sigma$ of $U(X)$ lies in $1+\fkm_A$, because we have $[\varpi^m](X)\equiv P_\phi(X)\equiv X^{q^{mn}}\ (\bmod\ \fkm_A)$. 

Also, once we fix a uniformizer $\varpi$ of $\mO$, we have a standard basis $\{e_1,\ldots,e_n\}$ of $(\fkp^{-m}/\mO)^n$ as a free $\mO/\fkp^m$-module, where $e_i:=(0,\ldots,0,\varpi^{-m},0,\ldots,0)$. We call $\phi(e_1),\ldots,\phi(e_n)\in \fkm_A$ the {\em formal parameters} of $\phi$. 

For $A\in \mC$, a {\em deformation of $\Sigma_n$ with level $\fkp^m$ structure} over $A$ is a triple $(\Sigma,i,\phi)$ consisting of a usual deformation $(\Sigma,i)$ over $A$ and a Drinfeld level $\fkp^m$ structure $\phi$ of $\Sigma$. Two deformations $(\Sigma,i,\phi), (\Sigma',i',\phi')$ are {\em equivalent} if there is an equivalence $f$ of $(\Sigma,i)$ to $(\Sigma',i')$ such that $\phi'=f_\fkm\circ \phi$. The deformation functor $\mF_m$ from $\mC$ to the category of sets is defined by sending $A\in \mC$ to the set of equivalence classes of deformations of $\Sigma_n$ with level $\fkp^m$ structure over $A$.

\bpr \label{deflev} \label{defm}
{\rm (Drinfeld \cite{Dr})}\ For every integer $m\geq 1$, the deformation functor $\mF_m$ is represented by an $n$-dimensional regular local ring $A_m$. The local $W$-algebra homomorphism $A_0\ra A_m$, representing the obvious forgetting morphism of functors $\mF_m\ra \mF_0$, is finite and flat, and the universal object over $A_m$ is a level $\fkp^m$ structure $\phi_m$ on $\Sigma^\univ_n\otimes_{A_0}A_m$. The formal parameters $X_1,\ldots,X_n\in \fkm_{A_m}$ of the universal level $\fkp^m$ structure $\phi_m$ gives a set of regular parameters of $A_m$.
\epr

This functor $\mF_m$ naturally factors through the category of sets with right $GL_n(\mO/\fkp^m)$-action as follows. For $A\in \mC$, the set $\mF_m(A)$ has a natural right $GL_n(\mO/\fkp^m)$-action induced by the action of $g\in GL_n(\mO/\fkp^m)$ on the deformations defined as
\[ (\Sigma,i,\phi)\longmapsto (\Sigma,i,\phi\circ g).\]
This results in a right $GL_n(\mO/\fkp^m)$-action on $\Spec A_m$. In terms of formal parameters $X_1,\ldots,X_n$ of the universal level $\fkp^m$ structure on $A_m$, this coincides with the left action given by the ``linear" action of $GL_n(\mO/\fkp^m)$ on the row vector $(X_1,\ldots,X_n)$ from the right, as elements of $\fkm_{\Sigma^\univ_n\otimes A_m}$, i.e.\ the addition and the $\mO/\fkp^m$-multiplications on $X_i$'s are the operations of $\Sigma^\univ_n\otimes A_m$. The finite flat covering $\Spec A_m\ra \Spec A_0$ is a Galois etale covering on the generic fibers, with the Galois group $GL_n(\mO/\fkp^m)$.


\subsection{Realization as a complete local ring of a Shimura variety}

Here we briefly recall from \cite{HT} the realization of this deformation ring as a complete local ring of certain Shimura variety, especially the proof of Lemma II.2.7, given in p.114 of the book. The detailed definition of the relevant Shimura variety is described carefully in \cite{HT}, which we omit here as it has little importance for us. In \cite{HT} Section III.4, the proper flat integral model $X_{U^p,m}$ of some Shimura variety over $\Spec {\mathcal O}_{F,w}$ is defined, where $U^p$ is a compact open subgroup of $G(\A^{\infty,p})$ for certain reductive group $G$ over $\Q$, and $m=(m_1,\ldots,m_r)\in (\Z_{\geq 0})^r$ is a multi-index. 

Now we can take the local field $F_w$ to be our $K$, and the index $m_1$ to be our $m$. Writing $\Sh:=X_{U^p,m}$, Lemma III.4.1(1) of \cite{HT} tells us that the completion $\hat{\mO}_{\Sh,\ol{s}}$ of the strict local ring of $\Sh$ at any geometric point $\ol{s}$ centered at a closed point $s$ with $h(s)=0$ (a ``supersingular" point, i.e.\ where the etale height $h(s)$ of the corresponding Barsotti-Tate group is 0), is isomorphic to our universal deformation ring $A_m$. Note that the existence of such a closed point is ensured by Lemma III 4.3 of \cite{HT}. This will be used in Section 4.2.



\section{The level $\fkp$ deformation space and the first blow-up}

We fix $n\geq 1$, and we are interested in the deformation space $X:=\Spec A_1$ of formal $\mO$-modules over $\ol{k}$ of height $n$ with level $\fkp$ structures. This $X$ is a regular flat scheme over $S:=\Spec W$ of relative dimension $n-1$ with a (formally) smooth generic fiber. We denote the universal formal $\mO$-module over $A_0$ by $\Sigma^\univ:=\Sigma^\univ_n$. From this section on, we denote the ring $A_1$ simply by $A$, and the maximal ideal of $A$ by $\fkm:=(X_1,\ldots,X_n)$, where $X_i$ are the formal parameters of the universal level $\fkp$ structure on $\Sigma^\univ\otimes A$, associated to our choice of uniformizer $\varpi$ of $\mO$.

\subsection{The equation of the space}

We start by computing the defining equation of this space $X=\Spec A$. By Propositoin \ref{deflev}, the formal parameters $X_1,\ldots,X_n$ give a set of regular parameters, therefore we have a surjective local homomorphism of local $W$-algebras
\[\tilde{A}:=W[[\widetilde{X}_1,\ldots,\widetilde{X}_n]]\lra A,\ \ \ \widetilde{X}_i\longmapsto X_i\ (1\leq i\leq n),\]
where $(\tilde{A},\tilde{\fkm})\in \mC$ is a formal power series ring in $\tilde{X}_1,\ldots,\tilde{X}_n$ over $W$ with the maximal ideal $\widetilde{\fkm}:=(\varpi,\widetilde{X}_1,\ldots,\widetilde{X}_n)$.

Let $I$ be the kernel of this surjection, so that $A\cong \tilde{A}/I$. Note that $\tilde{A}$ is an $(n+1)$-dimensional regular local ring, and $A$ is regular by Proposition \ref{deflev}. Hence the ideal $I$ has height one, therefore principal, generated by an element $t$ which is part of a system of regular parameters of $\tilde{A}$, i.e.\ an element $t\in I\setminus \widetilde{\fkm}^2$ (\cite{Mat}, Theorem 14.2). But for any element $t'\in I\setminus \widetilde{\fkm}^2$, we have $t'\widetilde{\fkm}=(t')\cap \widetilde{\fkm}^2$, thus the map $(t')/t'\widetilde{\fkm}\ra I/\widetilde{\fkm}I$ is an injection between 1-dimensional $\ol{k}$-vector spaces, therefore an isomorphism, which in turn gives
$(t')=I$ by Nakayama's lemma. Therefore, to determine $I$ we only need to find any element $t'\in I\setminus \widetilde{\fkm}^2$, and for this we observe:

\bpr \label{firsteq}
We have the following equality in $A$:
\[\varpi=u\cdot \!\!\!\prod_{\ul{a}\in k^n\setminus \{\ul{0}\}}\bigl([a_1](X_1)+_{\Sigma^\univ}\cdots +_{\Sigma^\univ}[a_n](X_n)\bigr)\]
with $u\in 1+\fkm\subset A^\times$, where $+_{\Sigma^\univ},[\cdot]$ denote the operations in $\fkm_{\Sigma^\univ\otimes A}$, and we used the notation $\ul{a}:=(a_1,\ldots,a_n)$ and $\ul{0}:=(0,\ldots,0)$. (Note that $[a_i](X_i)$ are well-defined because $[\varpi](X_i)=0$.) In particular, we have $\varpi\in \fkm^{q^n-1}$.
\epr

\bprf
For the universal Drinfeld structure $(\Sigma^\univ\otimes A,i,\phi)$ over $A$, by definition we have:
\[[\varpi](T)=U(T)P_\phi(T),\ \ U(T)\in A[[T]],\ \ \ \ P_\phi(T)=\!\!\!\prod_{x\in (\fkp^{-1}/\mO)^n}\!\!(T-\phi(x)).\]
As remarked before in Section 3.3, the constant term $u$ of $U(T)$ is in $1+\fkm$. By comparing the leading terms (i.e.\ the coefficients of $T$), we have the equality:
\[\varpi=u\cdot \prod_{x\in (\fkp^{-1}/\mO)^n\setminus \{\ul{0}\}}\phi(x)\]
in $A$. The definition of the formal parameters gives:
\begin{align*}
\prod_{x\in (\fkp^{-1}/\mO)^n\setminus \{\ul{0}\}}\phi(x) &= \prod_{\ul{a}\in k^n\setminus \{\ul{0}\}}\phi\bigl(a_1e_1+\cdots +a_ne_n\bigr)\\
&= \prod_{\ul{a}\in k^n\setminus \{\ul{0}\}}\bigl([a_1](X_1)+_{\Sigma^\univ}\cdots +_{\Sigma^\univ}[a_n](X_n)\bigr).
\end{align*}
\eprf

In order to define an element $t\in \tilde{A}$ which vanishes in $A$, first we lift the universal formal $\mO$-module $\Sigma^\univ$ to $\tilde{A}$. By definition of the deformation space, defining such a lift $\tilde{\Sigma}^\univ$ over $\tilde{A}$ is equivalent to giving a local $W$-algebra homomorphism $f$ making the following diagram commute:
\[\xymatrix{
\tilde{A} \ar[r] & A\\
 & A_0=W[[T_1,\ldots,T_{n-1}]] \ar@{.>}[ul]^f \ar[u]
}\]
where the right vertical arrow is the canonical map defined in Proposition \ref{defm}. The existence of such an $f$ is ensured by the formal smoothness of $A_0$ over $W$ (finding such $f$ only amounts to defining the images of $T_i$ by finding a power series of $X_1,\ldots,X_n$ with coefficients in $W$ which represent the images of $T_i$ in $A$). 

\bde
We choose and fix one such $f$, and define $\tilde{\Sigma}^\univ:= \Sigma^\univ\otimes_{A_0,f}\tilde{A}$.
\ede

Note that this does not mean that $\tilde{X}_1,\ldots,\tilde{X}_n$ are formal parameters of a level $\fkp$ structure on $\tilde{\Sigma}^\univ$ (that would amount to giving a section $A\ra \tilde{A}$), thus $[\varpi]_{\tilde{\Sigma}^\univ}(\tilde{X}_i)$ does not vanish and $[a](\tilde{X}_i)$ for $a\in k$ is not well-defined.

\bde \label{pai}
For each $a\in k^\times$, let $\tilde{a}\in \MU_{q-1}\subset \mO^\times$ be its multiplicative lift, and set $\tilde{0}=0$. For each $\ul{a}=(a_1,\ldots,a_n)\in k^n\setminus \{\ul{0}\}$, define a formal power series:
\[P_{\ul{a}}(\widetilde{X}_1,\ldots,\widetilde{X}_n):=[\tilde{a}_1](\widetilde{X}_1)+_{\widetilde{\Sigma}^\univ}\cdots +_{\widetilde{\Sigma}^\univ}[\tilde{a}_n](\widetilde{X}_n)\in \tilde{\fkm}\subset \tilde{A},\]
where $+_{\widetilde{\Sigma}^\univ},[\cdot]$ denotes the operations in $\widetilde{\fkm}_{\tilde{\Sigma}^\univ}$. 
\ede

We record some easy properties of these power series $P_{\ul{a}}$.

\bpr \label{paiprop}
\benu
\item $P_{\ul{a}}\equiv \tilde{a}_1\widetilde{X}_1+\cdots +\tilde{a}_n\widetilde{X}_n \pmod{(\widetilde{X}_1,\ldots,\widetilde{X}_n)^2}.$
\item If $\ul{a}=c\cdot \ul{a}'$ for $c\in k^\times$, then $P_{\ul{a}}=u_c\cdot P_{\ul{a}'}$ for a unit $u_c\in \tilde{A}^\times$.
\item If $a_{j+1}=\cdots =a_n=0$ for some $j<n$, then $P_{\ul{a}}\in (\widetilde{X}_1,\ldots,\widetilde{X}_j)$.
\eenu
\epr

\bprf
Parts (i) and (iii) follow from the definition of formal $\mO$-modules and $\tilde{0}=0$. For (ii), note that $P_{\ul{a}}=[\tilde{c}](P_{\ul{a}'})$, but as $\tilde{c}\in \mO^\times$, we have $[\tilde{c}](x)=u_c\cdot x$ for any $x\in \widetilde{\fkm}$, with a unit $u_c\in \tilde{c}+\widetilde{\fkm}$ depending on $x$.
\eprf

Now we can define an element of $\tilde{A}$ reducing to 0 in $A$ using Proposition \ref{firsteq}:

\bpr \label{localeq}
We have an $W$-algebra isomorphism
\[\tilde{A}/(P-\varpi)\cong A,\]
where $\tilde{A}:=W[[\widetilde{X}_1,\ldots,\widetilde{X}_n]]$, and $P\in \tilde{A}$ is a formal power series of the form:
\[P(\widetilde{X}_1,\ldots,\widetilde{X}_n) := \tilde{u}\cdot \!\!\!\prod_{\ul{a}\in k^n\setminus \{\ul{0}\}}P_{\ul{a}}(\widetilde{X}_1,\ldots,\widetilde{X}_n)\]
with $\tilde{u}\in 1+\widetilde{\fkm}\subset \tilde{A}^\times$.
\epr

\bprf
We define $\tilde{u}$ as any lift of $u\in A$ in Proposition \ref{firsteq} to $\tilde{A}$, which is a unit satisfying the asserted property because $\tilde{A}\ra A$ is a local homomorphism. The element $P-\varpi\in \tilde{A}$ reduces to 0 in $A$ by Proposition \ref{firsteq}, and by the discussion in the beginning of this section, we only need to make sure that $P-\varpi\notin \widetilde{\fkm}^2=(\varpi,\widetilde{X}_1,\ldots,\widetilde{X}_n)^2$. This is clear because $P\in (\widetilde{X}_1,\ldots,\widetilde{X}_n)^{q^n-1}$ and $P-\varpi\equiv -\varpi\ \pmod{\widetilde{\fkm}^2}$ unless when $q=2$ and $n=1$, in which case $P-\varpi\equiv \tilde{u}\widetilde{X}_1-\varpi$.
\eprf

\brem
{\rm This argument can also be applied to the deformation spaces for higher level $\fkp^m$ structures, giving similar equations with $P$ of the form
\[P=(\text{unit})\cdot \!\!\!\prod_{\ul{a}\in (\mO/\fkp^m)^n\setminus(\fkp/\fkp^m)^n}P_{\ul{a}}.\]
}
\erem

Here we recall from Section 2.3 that the left action of $(a_{ij})\in GL_n(k)$ on $A$ on the formal parameters $X_1,\ldots,X_n$ is defined by
\begin{equation} \label{glactionxi}
X_j\longmapsto [a_{1j}](X_1)+_{\Sigma^\univ}\cdots +_{\Sigma^\univ}[a_{nj}](X_n),
\end{equation}
therefore by definition we can regard it as:
\begin{equation} \label{glaction}
\widetilde{X}_j \bmod I\longmapsto P_{(a_{ij})_i}(\widetilde{X}_1,\ldots,\widetilde{X}_n) \bmod I
\end{equation}
where $(a_{ij})_i:=(a_{1j},\ldots,a_{nj})$ and $I:=(P-\varpi)$. Note that this does not lift to the action of $GL_n(k)$ on $\tilde{A}$, as $\tilde{a}$ are not additive lifts.

\subsection{The special fiber}

Now we investigate the special fiber of $X=\Spec A$. As $X$ is a scheme over $S:=\Spec W$, we use the notation $s:=\Spec \ol{k}$ and $X_s:=X\times_Ss$. By Proposition \ref{localeq}, we have:
\begin{equation} \label{spfibereq}
X_s=\Spec\ol{k}[[\widetilde{X}_1,\ldots,\widetilde{X}_n]]/\bigl(\prod_{\ul{a}\in k^n\setminus\{\ul{0}\}}(P_{\ul{a}}\bmod \fkp)\bigr).
\end{equation}

\bde
For each ${\ul{a}}\in k^n\setminus\{\ul{0}\}$, we denote by $Y_{\ul{a}}$ the closed subscheme of $X_s$ defined by $(P_{\ul{a}} \bmod \fkp)=0$, or equivalently, the closed subscheme of $X$ defined by $(P_{\ul{a}} \bmod I)=0$ (note that $(P_{\ul{a}}\bmod I)$ divides $\varpi$ in $A$).
\ede

Note that, in the regular local ring $\ol{k}[[\tilde{X}_1,\ldots,\tilde{X}_n]]$ with the maximal ideal $\ol{\fkm}:=(\tilde{X}_1,\ldots,\tilde{X}_n)$, we have
\[ (P_{\ul{a}}\bmod \fkp) \equiv a_1\tilde{X}_1 +\cdots + a_n\tilde{X}_n\ (\bmod\ \ol{\fkm}^2) \]
by Proposition \ref{paiprop}(i), therefore that $(P_{\ul{a}}\bmod \fkp) \in \ol{\fkm}\setminus \ol{\fkm}^2$. This shows that the quotient $\ol{k}[[\tilde{X}_1,\ldots,\tilde{X}_n]]/(P_{\ul{a}}\bmod \fkp)$ is an $(n-1)$-dimensional regular local ring, thus its spectrum $Y_{\ul{a}}$ is an irreducible and reduced $\ol{k}$-scheme of dimension $n-1$.

By Proposition \ref{paiprop}(ii), we see that $Y_{\ul{a}}=Y_{\ul{a}'}$ when $\ul{a}=c\cdot \ul{a}'$ for $c\in k^\times$, and Proposition \ref{paiprop}(i) ensures that $Y_{\ul{a}}\neq Y_{\ul{a}'}$ if otherwise. Therefore we introduce the following notation to label these closed subschemes of $X_s$, by $k$-rational hyperplanes of a projective $(n-1)$-space. 

\bde \label{Pslabel}
Let $\Ps:=\Ps^{n-1}$ be an $(n-1)$-dimensional projective space over $\ol{k}$, equipped with a set of projective coordinates $(X_1^*:\cdots :X^*_n)$. For each ${\ul{a}}\in k^n\setminus\{\ul{0}\}$, we define a $k$-rational hyperplane:
\[M_{\ul{a}}:a_1X^*_1+\cdots +a_nX^*_n=0\]
of $\Ps$. For any such hyperplane $M=M_{\ul{a}}$, we denote $Y_M:=Y_{\ul{a}}$, which is well-defined by the above remark.
\ede

We will later identify this $\Ps$ with the exceptional divisor of our first blow-up. The following proposition is clearly seen from the equation (\ref{spfibereq}) and the remarks before Definition \ref{Pslabel}.

\bpr
The correspondence $M\mapsto Y_M$ gives a bijection from the set of $k$-rational hyperplanes of $\Ps$ to the set of irreducible components of $X_s$. In particular, there are $(q^n-1)/(q-1)$ irreducible components of $X_s$. Each irreducible component $Y_M$ has multiplicity $|k^\times|=q-1$ in $X_s$.
\epr

We will also need the closed subschemes of $X_s$ of higher codimension, i.e.\ the intersections of $Y_M$'s, in Section 4.

\bde
\benu
\item For any $k$-rational linear subspace $N\subset \Ps$, we define a reduced closed subscheme $Y_N$ of $X_s$ as 
\[Y_N:=\bigcap_{N\subset M}Y_M.\]
\item For $1\leq h\leq n-1$, let $Y^{[h]}:=\bigcup_NY_N$, where $N$ runs through all the $k$-rational linear subspaces of $\Ps$ of dimension $h-1$. Define $Y^{[0]}$ as the unique closed point $x$ of $X$. Note that $Y^{[n-1]}=\bigcup_MY_M=X_s^{\rm red}$.
\eenu
\ede

The way closed subschemes $Y_M$ intersect with each other correspond precisely to how the hyperlanes $M$ intersect inside $\Ps$. More precisely, we have the following.

\blem \label{yn}
Let $1\leq h\leq n-1$, and $N\subset \Ps$ be a $k$-rational linear subspace of dimension $h-1$. If we choose a set of linearly independent $k$-rational hyperplanes $M_1,\ldots,M_{n-h}$ such that $N=\bigcap_{i=1}^{n-h}M_i$ in $\Ps$, then we have $Y_N=\bigcap_{i=1}^{n-h}Y_{M_i}$. 

Therefore the defining ideal of $Y_N$ is $(P_{\ul{a}_1},\ldots,P_{\ul{a}_{n-h}})$ if $M_i=M_{\ul{a}_i}$. In particular, the closed immersion $Y_N\subset X$ is a regular immersion of codimension $n-h$, hence $Y_N$ is a smooth $\ol{k}$-scheme of dimension $h$.
\elem

\bprf
As the other inclusion is trivial, it is enough to show $Y_N\supset \bigcap_{i=1}^{n-h}Y_{M_i}$, or $Y_M\supset \bigcap_{i=1}^{n-h}Y_{M_i}$ for any $k$-rational hyperplane $M$ containing $N$. For this, by the $GL_n(k)$-action we can assume that $N\subset \Ps$ is defined by $X_1^*=\cdots=X_{n-h}^*=0$ and $M_i$ is defined by $X_i^*=0$ for $1\leq i\leq n-h$, without loss of generality. In this case $Y_{M_i}$ is defined by $P_{\ul{1_i}}=\widetilde{X}_i=0$ where $\ul{1_i}=(0,\ldots,0,1,0,\ldots,0)$ with 1 in the $i$-th entry, therefore $\bigcap_{i=1}^{n-h}Y_{M_i}$ is defined by the ideal $(X_1,\ldots,X_{n-h})$ of $A$. Now for any $M$ containing $N$ is defined by an equation $\sum_{i=1}^na_iX_i^*=0$ with $a_{n-h+1}=\cdots =a_n=0$, which shows $P_{\ul{a}}\in (\widetilde{X}_1,\ldots,\widetilde{X}_{n-h})$ by Proposition \ref{paiprop}(iii), hence $Y_M\supset \bigcap_{i=1}^{n-h}Y_{M_i}$.
\eprf

The action (\ref{glaction}) of $GL_n(k)$ on the formal parameters $X_1,\ldots,X_n$ and the Definition \ref{pai} of $P_{\ul{a}}$ show that the right $GL_n(k)$-action on $X$ permutes the components $Y_M$ in an obvious way, by acting on the set of indices through $PGL_n(k)$. More precisely, for $g\in GL_n(k)$ and $\ul{a}\in k^n\setminus\{\ul{0}\}$, let $\ul{a}\mapsto \ul{a}g$ be the linear action on the row vector $\{\ul{a}\}$. This gives the action $M\mapsto Mg$ on the hyperplanes of $\Ps$, induced by the right linear $PGL_n(k)$-action on $\Ps$ through the projective coordinates. Then $GL_n(k)$ acts on the set of irreducible components of $X_s$ by $Y_M\mapsto Y_{Mg}$. Similarly we see that $Y^{[h]}$ is stable under this action which maps $Y_N\mapsto Y_{Ng}$.

\subsection{Some preliminaries on blow-ups of schemes}

Before we start blowing up our space, we collect some facts on blow-up of schemes which will be used in the sequel. Let $X$ be a noetherian scheme and $\mI$ be a coherent sheaf of ideals in $\mO_X$, and $Y:=\Spec(\mO_X/\mI)$ be the closed subscheme defined by $\mI$. The {\it blow-up} of $X$ at $Y$ is defined as the scheme $X':=\Proj\bigl(\bigoplus_{i\geq 0}\mI^i\bigr)$, which is projective over $X$. The structure morphism $p:X'\ra X$ is an isomorphism outside the inverse image $Y':=p^{-1}(Y)$ of $Y$ (the {\it exceptional divisor}). As a projective scheme over $Y$, we have $Y'=\Proj\bigl(\bigoplus_{i\geq 0}\mI^n/\mI^{n+1}\bigr)$. Therefore, if $Y\subset X$ is a {\it regular immersion} of codimension $r$, in other words when $\mI$ is locally generated by a regular sequence of sections of $\mO_X$ of length $r$, then the graded $\mO_X/\mI$-algebra $\bigoplus_{i\geq 0}\mI^n/\mI^{n+1}$ is naturally isomorphic to the symmetric algebra generated by the conormal sheaf $\mN_{Y/X}=\mI/\mI^2$, which is a locally free $\mO_X/\mI$-module of rank $r$, thus $Y'\cong \Ps(\mN_{Y/X})$ is a projective space bundle ($\Ps^{r-1}$-bundle) over $Y$ (\cite{EGA} IV, 16.9 and 19.4), and if $\mI$ is locally generated by a regular sequence $X_1,\ldots,X_r$, then $(X^*_1:\cdots:X^*_r)$ gives a set of projective coordinates of this projective space bundle where $X^*_i:=X_i\bmod\mI^2$.

We will need some commutativity between blow-ups and completions.

\blem \label{blbc}
Let $X,\mI,X'$ as above, and $f:Z\ra X$ be a flat morphism. Let $\mI':=\im(f^*\mI\ra \mO_Z)$, where $f^*\mI=f^{-1}\mI\otimes_{f^{-1}\mO_X}\mO_Z$, and $Z'$ be the blow-up of $Z$ at the subscheme defined by $\mI'$. Then the canonical morphism $Z'\ra Z\times_XX'$ defined by the universality of blow-ups is an isomorphism.
\elem

\bprf
As $f$ is flat, the sheaf $\mO_Z$ is a flat $f^{-1}\mO_X$-algebra, which shows that $f^*\mI\ra \mO_Z$ is an injection and $f^*\mI\cong \mI'$. Now $Z\times_XX'=Z\times_X\Proj\bigl(\bigoplus_{i\geq 0}\mI^i\bigr)=\Proj\bigl(\mO_Z\otimes_{f^{-1}\mO_X}\bigoplus_{i\geq 0}(f^{-1}\mI)^i\bigr)=\Proj\bigl(\bigoplus_{i\geq 0}(f^*\mI)^i\bigr)$, which is isomorphic to $Z'=\Proj\bigl(\bigoplus_{i\geq 0}(\mI')^i\bigr)$.
\eprf

\bcor \label{blcp}
Let $X,\mI$ as above. For $x\in X$, let $\widehat{X}_x:=\Spec \widehat{\mO}_{X,x}$ be the complete local ring (or the completion of the strict local ring) of $X$ at $x$, and $X'_x$ the blow-up of $\widehat{X}_x$ at the subscheme defined by $\mI\cdot \mO_{\widehat{X}_x}$. Then $X'_x\cong \widehat{X}_x\times_XX'$.
\ecor

\subsection{The first blow-up $Z_1$}

We first blow-up the unique closed point $x=Y^{[0]}$ of $X$. By Proposition \ref{localeq}, we have the description 
\begin{gather*}
X=\Spec W[[\widetilde{X}_1,\ldots,\widetilde{X}_n]]/(P-\varpi)\\
P(\widetilde{X}_1,\ldots,\widetilde{X}_n)=\tilde{u}\cdot \!\!\!\prod_{\ul{a}\in k^n\setminus \{\ul{0}\}}P_{\ul{a}}(\widetilde{X}_1,\ldots,\widetilde{X}_n)
\end{gather*}
with $\tilde{u}\in 1+\widetilde{\fkm}\subset W[[\widetilde{X}_1,\ldots,\widetilde{X}_n]]^\times$. The universal formal parameters on $X$ are given by $X_i=\tilde{X_i}\bmod (P-\varpi)$ for $1\leq i\leq n$.
 
\bde \label{z0de}
Let $Z_1\ra X$ be the blow-up of $X$ at the unique closed point $x:=(X_1,\ldots,X_n)$ of $X$, and $Y_\Ps$ be its exceptional divisor, i.e.\ the inverse image of $x$ in $Z_1$.
\ede

As $X$ is the spectrum of an $n$-dimensional regular local ring, the exceptional divisor $Y_\Ps$ is an $(n-1)$-dimensional projective space $\Ps^{n-1}$ over $\ol{k}$, equipped with a set of projective coordinates $(X^*_1:\cdots:X^*_n)$, where $X_i^*:=X_i\bmod\fkm^2$. We identify $Y_\Ps$ with the $\Ps$ that we introduced in Definition \ref{Pslabel} for the indexing purpose. 

As the center $x$ of the blow-up is $GL_n(k)$-invariant, the blow-up $Z_1$ inherits the right $GL_n(k)$-action on $X$. The action on the exceptional divisor $Y_\Ps\cong \Ps$ is the obvious one:

\bpr \label{actz0}
The right action of $GL_n(k)$ on $Y_\Ps\cong \Ps$, in terms of the projective coordinates $(X^*_1:\dots:X^*_n)$, is the right linear action through $PGL_n(k)$.
\epr

\bprf
As the action of $GL_n(k)$ on $Y_\Ps\cong \Ps(\fkm/\fkm^2)$ is induced from the action on the $n$-dimensional $k$-vector space $\fkm/\fkm^2$ generated by $X_i^*=X_i \bmod \fkm^2$ for $1\leq i\leq n$, it suffices to see that the left action of $(a_{ij})\in GL_n(k)$ on $X_i^*$ is defined by
\[X_j^*\longmapsto a_{1j}X_1^*+\cdots +a_{nj}X_n^*,\]
which is readily confirmed by reducing the left action (\ref{glactionxi}) on $A$ modulo $\fkm^2$.
\eprf

Now we analyze how the proper transforms of the other components $Y_M$ intersect the exceptional divisor $Y_\Ps$ in $Z_1$. 

\bde
For a $k$-linear subspace $N\subset \Ps$ of dimension $h-1$, with $1\leq h\leq n-1$, we denote the proper transform of $Y_N$ in $Z_1$ by $Y_{N,1}$. 

For $1\leq h\leq n-1$, let $Y^{[h]}_1$ be the proper transform of $Y^{[h]}$ in $Z_1$. Then $Y^{[h]}_1=\bigcup_NY_{N,1}$ for $1\leq h\leq n-1$, where $N$ runs through all $k$-linear subspaces of $\Ps$ of dimension $h-1$.
\ede

\bpr \label{z0eq}
Let $Z_1$ be the model of $X$ defined in Definition \ref{z0de}.
\benu
\item In the special fiber of $Z_1$, the exceptional divisor $Y_\Ps$ has multiplicity $q^n-1$.
\item For a $k$-linear subspace $N\subset \Ps\cong Y_\Ps$, we have $Y_\Ps\cap Y_{N,1}=N$.
\eenu
\epr

\bprf
These properties can be checked by looking at the completions along the exceptional divisor of the affine open sets of $Z_1$. We can reduce to the case where $N$ is a hyperplane $M=M_{\ul{a}}$ in (ii), and assume moreover that $a_n\neq 0$, without loss of generality. We look at the affine open set that is the spectrum of  
\[W[[\widetilde{X}_1,\ldots,\widetilde{X}_n]][V_1,\ldots,V_{n-1}]/(P(\widetilde{X}_1,\ldots,\widetilde{X}_n)-\varpi,\ V_i\widetilde{X}_n-\widetilde{X}_i).\]
We take the completion of this ring along the exceptional divisor $\widetilde{X}_n=0$ (which is an affine space $\A^{n-1}$ with the coordinate ring $k[V_1,\ldots,V_{n-1}]$) to get the spectrum of
\[B_1:=W[V_1,\ldots,V_{n-1}][[\widetilde{X}_n]]/(P(V_1\widetilde{X}_n,\ldots, V_{n-1}\widetilde{X}_n, \widetilde{X}_n)-\varpi).\]
with $P=\tilde{u}\cdot\prod_{\ul{a}\in k^n\setminus\{\ul{0}\}}P_{\ul{a}}$ where
\[\tilde{u}\in W[V_1,\ldots,V_{n-1}][[\widetilde{X}_n]]^\times,\ \ \tilde{u}\equiv 1\ (\bmod\ \widetilde{X}_n).\]
Now by the congruence Proposition \ref{paiprop}(i), we have
\[P_{\ul{a}}(V_1\widetilde{X}_n,\ldots,V_{n-1}\widetilde{X}_n,\widetilde{X}_n) = \widetilde{X}_n\cdot P'_{\ul{a}}\]
with some $P'_{\ul{a}}\in W[V_1,\ldots,V_{n-1}][[\widetilde{X}_n]]$, which satisfies
\begin{equation} \label{z0eq1}
P'_{\ul{a}}(V_1,\ldots,V_{n-1},\widetilde{X}_n) \equiv \sum_{i=1}^{n-1}\tilde{a_i}V_i+\tilde{a_n} \pmod{\widetilde{X}_n},
\end{equation}
hence we can write
\[B_1=W[V_1,\ldots,V_{n-1}][[\widetilde{X}_n]]/\bigl(\tilde{u}\cdot \widetilde{X}_n^{q^n-1}\cdot \prod_{\ul{a}\in k^n\setminus\{\ul{0}\}}P'_{\ul{a}}-\varpi\bigr).\]
Therefore the special fiber of this affine scheme is the spectrum of:
\[B_1\otimes_W\ol{k} = \ol{k}[V_1,\ldots,V_{n-1}][[\widetilde{X}_n]]/\bigl(\widetilde{X}_n^{q^n-1}\cdot \prod_{\ul{a}\in k^n\setminus\{\ul{0}\}}(P'_{\ul{a}}\bmod \fkp)\bigr).\]
In this expression, the exceptional divisor $\widetilde{X}_n=0$ clearly has multiplicity $q^n-1$ in the special fiber $\Spec (B_1\otimes_W\ol{k})$. Also the inverse image in $\Spec B_1$ of the proper transform $Y_{M,1}$ of $Y_M$ for $M=M_{\ul{a}}$ is defined by $P'_{\ul{a}}(V_1,\ldots,V_{n-1},\widetilde{X}_n)=0$, therefore by the above congruence it intersects $\widetilde{X}_n=0$ at the hyperplane $\sum_{i=0}^{n-1}a_iV_i+a_n=0$ in $\A^{n-1}$. As the situations in the completions of the other affine open sets of $Z_1$ are checked in exactly the same way, we conclude that $Y_{M,1}$ intersects $Y_\Ps$ at the hyperplane $M_{\ul{a}}\subset \Ps\cong Y_\Ps$.
\eprf

The following simple stratification $Y_\Ps^{(h)}$ on $Y_\Ps\cong \Ps$ by locally closed subschemes of dimension $h-1$ will be used in the analysis of this model in the subsequent sections.

\bde
For $1\leq h\leq n$, let $Y_\Ps^{[h]}$ denote the reduced closed subscheme of $Y_\Ps\cong \Ps$ which is the union of all $k$-rational linear subspaces of dimension $h-1$. Then by Proposition \ref{z0eq}(ii), we have $Y_\Ps\cap Y^{[h]}_1=Y_\Ps^{[h]}$ for $1\leq h\leq n-1$. 

Also, set $Y_\Ps^{[0]}:=\emptyset$, and $Y_\Ps^{(h)}:=Y_\Ps^{[h]}-Y_\Ps^{[h-1]}$ for $1\leq h\leq n$.
\ede

Then $Y_\Ps^{(h)}$ is a smooth locally closed subscheme of $Y_\Ps\cong \Ps$ of dimension $h-1$, namely the disjoint union
\[Y_\Ps^{(h)}=\coprod_NN^0,\]
where $N\subset \Ps$ runs through the set of $k$-rational linear subspaces of dimension $h-1$, and $N^0$ is the complement in $N$ of all $k$-rational proper linear subspaces of $N$. This stratification is clearly stable under the $GL_n(k)$-action. 

\section{A generalized semistable model $Z_\st$}

Now we construct a generalized semistable model of $X$ over $S$ by blowing up the model $Z_1$ further. Recall that, in the special fiber of $Z_1$, there is one proper component $Y_\Ps\cong \Ps$ with multiplicity $q^n-1$ and a set of projective coordinates $(X^*_1:\cdots:X^*_n)$ on it, and $(q^n-1)/(q-1)$ other components, each with multiplicity $q-1$, intersecting $Y_\Ps$ at each of the $k$-rational hyperplanes. While the intersections of reduced components are pairwise transversal, more than $n$ components meet at a point, so we need to blow-up several more times to get a generalized semistable model. The successive blow-up performed here is analogous to the one performed to obtain the irreducible components of the special fiber of the $p$-adic upper half space (see for example \cite{It}, Section 4). 

\subsection{A generalized semistable model}

Recall that $Z_1$ was the blow-up of $X$ at the closed point $x=Y^{[0]}$. 

\bde
We define inductively $Z_{h+1}\ra Z_h$ for $1\leq h\leq n-2$ as the blow-up of $Z_h$ at the proper transform of $Y^{[h]}$ in $Z_h$. We write $Z_\st:=Z_{n-1}$.
\ede

The scheme $Z_2$ is the blow-up of $Z_1$ at $Y^{[1]}_1=\coprod_NY_{N,1}$, where $N$ runs through all the $k$-rational points of $\Ps$. Similarly, the blow-up $Z_{h+1}\ra Z_h$ is centered at $\coprod_NY_{N,h}$, where $N$ runs through all the $k$-rational linear subspaces of $\Ps$ of dimension $h-1$ and $Y_{N,h}$ denotes the proper transform of $Y_N$ in $Z_h$. 

For a $k$-rational point $N$, as $Y_N$ has codimension $n-1$ in $Z_1$, the inverse image $Y_{N,2}$ of $Y_{N,1}$ in $Z_2$ is a $\Ps^{n-2}$-bundle over $Y_{N,1}$, which is an irreducible component of the special fiber of $Z_2$. As $(q^{n-1}-1)/(q-1)=|\Ps^{n-2}(k)|$ components meet at $Y_N$, the component $Y_{N,2}$ has multiplicity $q^{n-1}-1$ in the special fiber of $Z_2$. Similarly, for a $k$-rational linear subspace $N$ of dimension $h-1$, with $1\leq h\leq n-1$, the proper transform $Y_{N,h}$ in $Z_h$ is a regular immersion of codimension $n-h$ in $Z_h$, hence its inverse image $Y_{N,h+1}$ in $Z_{h+1}$ is a $\Ps^{n-h-1}$-bundle over $Y_{N,h}$. 

Eventually, for a $k$-rational linear subspace $N$ of dimension $h-1$, with $1\leq h\leq n-1$, the proper transform $Y_{N,\st}$ of $Y_N$ in $Z_\st$ is an irreducible component of the special fiber of $Z_\st$, whose multiplicity is $q^{n-h}-1$. The inverse image $Y_{\Ps,\st}$ of $Y_\Ps$ in $Z_\st$ is the unique proper component of the special fiber of $Z_\st$, whose multiplicity is $q^n-1$. This $Y_{\Ps,\st}$ is isomorphic to the $B^{n-1}$ in the notation of \cite{It}, Section 4, and $Y_{N,\st}\cap Y_{\Ps,\st}$ is isomorphic to $B^h\times B^{n-h-2}$.

Note that all the blow-ups performed here are $GL_n(k)$-equivariant, so that the models $Z_1,Z_2,\ldots,Z_{n-1}=Z_\st$ all inherit the right $GL_n(k)$-action on $X$, and $g\in GL_n(k)$ sends $Y_{N,\st}$ to $Y_{Ng,\st}$.

Now our main result here is:

\bth \label{gsstmodel}
The $W$-scheme $Z_\st=Z_{n-1}$ is generalized semistable. This means that, at every closed point of $Z_\st$, its complete local ring is isomorphic to 
\[ W[[T_1,\ldots,T_n]]/(T_1^{e_1}\cdots T_d^{e_d}-\varpi)\ \ (d\leq n),\] 
where integers $e_i$ are all prime to $\chara k$.
\ethm

\bprf
We start by describing the complete local rings at all the closed points of $Z_1$, i.e.\ all the closed points of $Y_\Ps$. By the obvious $GL_n(k)$-symmetry, it is enough to look at each points on the completion along the exceptional divisor of the affine piece that was defined in the proof of Proposition \ref{z0eq}:
\[\Spec B_1 = \Spec W[V_1,\ldots,V_{n-1}][[\widetilde{X}_n]]\bigm/\bigl(\tilde{u}\cdot \widetilde{X}_n^{q^n-1}\cdot \!\!\!\prod_{\ul{a}\in k^n\setminus\{\ul{0}\}}P'_{\ul{a}}-\varpi\bigr), \]
where $\tilde{u}\in W[V_1,\ldots,V_{n-1}][[\widetilde{X}_n]]^\times$. Hereafter in this proof, we denote the coordinates of a point in the corresponding lower case alphabets to distinguish them from the elements of the coordinate rings. Let $x\in \Spec B_1$ be a point which is closed in $Z_1$ and $(v_1,\ldots,v_{n-1})\in \ol{k}^{n-1}$ the affine coordinates of $x$. The complete local ring $\widehat{\mO}_x$ depends on whether the value of $P'_{\ul{a}}$ at $x$ is a unit or not for each $\ul{a}$, which in turn depends on how many of the $v_i$'s have $k$-rational linear relations among them, by the equations (\ref{z0eq1}). 

If the image of $x$ in $Z_1$ lies in $Y_\Ps^{(n)}$, i.e.\ if there is no $k$-rational relation of the form
\[ a_1v_1+\cdots +a_{n-1}v_{n-1}+a_n=0\ \ \ (a_1,\ldots,a_n\in k)\]
among $v_1,\ldots,v_{n-1}$, then we have $\widehat{\mO}_x\cong W[[\widetilde{V}_1,\ldots,\widetilde{V}_{n-1},\widetilde{X}_n]]/(\tilde{u}\cdot \widetilde{X}_n^{q^n-1}-\varpi)$ with $\tilde{u}\in W[[\widetilde{V}_1,\ldots,\widetilde{V}_{n-1},\widetilde{X}_n]]^\times$, where $\widetilde{V_i}$ is a translation of $V_i$ which vanishes at $x$. If $x\in Y_\Ps^{(m)}$ with $m\leq n-1$, i.e.\ there are exactly $m$ of the $k$-rational linear dependences between $v_1,\ldots,v_{n-1}$, then by $GL_n(k)$-symmetry, it suffices to treat the case where $v_1=\cdots =v_m=0$ and no $k$-rational linear relation between $v_{m+1},\ldots,v_{n-1}$, by using the $GL_n(k)$-action. In this case the value of $P'_{\ul{a}}$ at $x$ is a non-unit if and only if $a_{m+1}=\cdots = a_n=0$. Thus we have 
\[\widehat{\mO}_x\cong W[[\widetilde{V}_1,\ldots,\widetilde{V}_{n-1},\widetilde{X}_n]]\bigm/\bigl(\tilde{u}\cdot \widetilde{X}_n^{q^n-1}\cdot\!\!\!\prod_{\ul{a}\in k^m\setminus\{\ul{0}\}}P'_{(\ul{a},\ul{0})}-\varpi\bigr)\]
with $\tilde{u}\in W[[\widetilde{V}_1,\ldots,\widetilde{V}_{n-1},\widetilde{X}_n]]^\times$, where $\widetilde{V_i}$ is a translation of $V_i$ which vanishes at $x$ and $(\ul{a},\ul{0}):=(a_1,\ldots,a_m,0,\ldots,0)$.

By Proposition \ref{paiprop}(i),(iii) we have $P'_{(\ul{a},\ul{0})}(\widetilde{V}_1,\ldots,\widetilde{V}_{n-1},\widetilde{X}_n) \in (\widetilde{V}_1,\ldots,\widetilde{V}_m)$, and 
\begin{equation} \label{p'cong}
P'_{(\ul{a},\ul{0})}(\widetilde{V}_1,\ldots,\widetilde{V}_{n-1},\widetilde{X}_n) \equiv \sum_{i=1}^{m}\tilde{a_i}\widetilde{V_i} \pmod{(\widetilde{V}_1,\ldots,\widetilde{V}_m)^2}.
\end{equation}
As a special case where $m=1$, we see that 
\[\widehat{\mO}_x\cong W[[\widetilde{V}_1,\ldots,\widetilde{V}_{n-1},\widetilde{X}_n]]\bigm/\bigl(\tilde{u}\cdot \widetilde{X}_n^{q^n-1}\cdot\prod_{a\in k^\times}P'_{(a,\ul{0})}-\varpi\bigr)\]
is already generalized semistable, because by Proposition \ref{paiprop} we have
\[\prod_{a\in k^\times}P'_{(a,\ul{0})}=(\text{unit})\cdot \widetilde{V}_1^{q-1}.\]

Now by Corollary \ref{blcp}, we can analyze the situation after the following blow-ups by blowing up these complete local rings further. As we go through the procedure of the successive blow-ups, the first time we touch the point $x\in Y_\Ps^{(n-m)}$ satisfying $v_1=\cdots =v_m=0$ is when we blow-up $Y_{N,n-m}$, where $N$ is the codimension $m$ linear subspace $X^*_1=\cdots =X^*_m=0$ on which $x$ lies on. Then by the proof of Lemma \ref{yn} we see that $Y_{N,n-m}$ is defined by the ideal $(V_1,\ldots,V_m)$, which is a regular sequence in $\widehat{\mO}_x$. After the blowing up and completing along the new exceptional divisor, we obtain the affine pieces of the form:
\begin{gather*}
\Spec W[U_1,\ldots,U_{m-1}][[\widetilde{V}_m,\ldots,\widetilde{V}_{n-1},\widetilde{X}_n]]\bigm/ I,\\
I=\bigl(\tilde{u}\cdot \widetilde{X}_n^{q^n-1}\cdot\!\!\!\prod_{\ul{a}\in k^m\setminus\{\ul{0}\}}P'_{(\ul{a},\ul{0})}(U_1\widetilde{V}_m,\ldots,U_{m-1}\widetilde{V}_m,\widetilde{V}_m,\ldots,\widetilde{V}_{n-1},\widetilde{X}_n)-\varpi\bigr),
\end{gather*}
with a unit $\tilde{u}$, and by the congruence (\ref{p'cong}) we can pull out the $V_m$ as we did in (\ref{z0eq1}):
\[P'_{(\ul{a},\ul{0})}(U_1\widetilde{V}_m,\ldots,U_{m-1}\widetilde{V}_m,\widetilde{V}_m,\ldots,\widetilde{V}_{n-1},\widetilde{X}_n) = \widetilde{V}_m\cdot P''_{\ul{a}}\]
for some $P''_{\ul{a}}\in W[U_1,\ldots,U_{m-1}][[\widetilde{V}_m,\ldots,\widetilde{V}_{n-1},\widetilde{X}_n]]$ which satisfies
\[P''_{\ul{a}}(U_1,\ldots,U_{m-1},\widetilde{V}_m,\ldots,\widetilde{V}_{n-1},\widetilde{X}_n) \equiv \sum_{i=1}^{m-1}\tilde{a_i}U_i+\tilde{a_m} \pmod{\widetilde{V}_m},\]
hence we can write the above affine piece in the form
\[\Spec W[U_1,\ldots,U_{m-1}][[\widetilde{V}_m,\ldots,\widetilde{V}_{n-1},\widetilde{X}_n]] \bigm/ \bigl(\tilde{u}\cdot \widetilde{X}_n^{q^n-1}\cdot \widetilde{V}_m^{q^m-1} \cdot \!\!\!\prod_{\ul{a}\in k^m\setminus\{\ul{0}\}}P''_{\ul{a}}-\varpi\bigr).\]
Then if we look at the closed point $y$ on the exceptional divisor of this blow-up, which is a $\Ps^{m-1}$-bundle over $Y_{N,n-m}$ with a set of projective coordinates $(\widetilde{V}^*_1:\cdots :\widetilde{V}^*_m)$ or affine coordinates $(u_1,\ldots,u_{m-1})$ with $U_i=\widetilde{V}^*_i/\widetilde{V}^*_m$ on the fibers, the complete local ring $\widehat{\mO}_y$ will depend on how many $k$-rational linear relations there are among the fiber coordinates $u_1,\ldots,u_{m-1}$ of $y$. Assume that $y$ lies on a $k$-rational linear subspace of codimension $l\leq m-1$, say $U_1=\cdots =U_l=0$ and no $k$-rational linear relation among $u_{l+1},\ldots,u_{m-1}$, without loss of generality because of the action of the parabolic subgroup of $GL_n(k)$ preserving $N$. Then the complete local ring will have the form:
\[\widehat{\mO}_y\cong W[[\widetilde{U}_1,\ldots,\widetilde{U}_{m-1},\widetilde{V}_m,\ldots,\widetilde{V}_{n-1},\widetilde{X}_n]]\bigm/
\bigl(\tilde{u}\cdot \widetilde{X}_n^{q^n-1}\cdot \widetilde{V}_m^{q^m-1} \prod_{\ul{a}\in k^l\setminus\{\ul{0}\}}P''_{(\ul{a},\ul{0})}-\varpi\bigr),\]
with $\tilde{u}\in W[[\widetilde{U}_1,\ldots,\widetilde{U}_{m-1},\widetilde{V}_m,\ldots,\widetilde{V}_{n-1},\widetilde{X}_n]]^\times$, and it is generalized semistable if $l=0$ or $1$. As we repeat the blow-up process, we get a sequence 
\[n=n_0,\ m=n_1,\ l=n_2,\ n_3,\ldots\]
until there is no more blow-up to be done when $n_t=1$ for some integer $t$. This proves that $Z_{n-1}$ is generalized semistable at all closed points.
\eprf

\subsection{Relation with a generalized semistable model of Shimura varieties}

In order to apply the computations of nearby cycle sheaves for the varieties with generalized semistable reduction due to T.\ Saito \cite{Sa} (we will recall his results in Section 6.2), we need to compare what we have constructed with some scheme of finite type over $S=\Spec W$. For this, first we can approximate the scheme $X$ by a scheme of finite type over $S$, so that its completion at the origin would be isomorphic to $X$. Then by performing the successive blow-ups of this scheme by the corresponding closed subschemes and observing that it has generalized semistable reduction (i.e.\ etale locally etale over $\Spec W[T_1,\ldots,T_n]/(T_1^{e_1}\cdots T_d^{e_d}-\varpi)\ \ (d\leq n)$), we can show that the results concerning the nearby cycles of schemes of finite type are applicable to our situation. We need the following:

\bpr
There is a scheme $U$ of finite type over $S$ and a closed point $x\in U$ such that $X$ is isomorphic to the spectrum of the complete local ring $\hat{\mO}_{U,x}$ of $U$ at $x$. Moreover, there is a closed subscheme $V^{[h]}$ of $U$ for each $0\leq h\leq n-1$, and $Y^{[h]}\cong X\times_UV^{[h]}$ via the morphism $X\cong \Spec \hat{\mO}_{U,x}\ra U$.
\epr

This proposition can be proven by a purely local argument. For the existence of $U$ and $x$, we can refer to \cite{St3}, Theorem 2.3.1, where the proof attributed to L.\ Fargues, using Artin's approximation theorem \cite{Ar} and Faltings' theory of strict $\mO$-modules \cite{Fa2}, is given. As this $U$ comes with an approximation of strict $\mO$-module $\Sigma^\univ[\fkp]$, its reduced closed subschemes $V^{[h]}=\bigcup_MV_M$, where $M$ runs through all $k$-subspaces of $(\fkp^{-1}/\mO)^n$ of dimension $n-h$, can be given in the exactly the same way as in the proof of Lemma \ref{yhshh} below (based on Lemma 9 of \cite{Man}). (Or, once we have an approximation of the finite covering $X\ra \Spec A_0 = \Spec W[[T_1,\ldots,T_{n-1}]]$ by a finite covering of the strict henselization of $W[T_1,\ldots,T_{n-1}]$ at the origin, we can use the fact that we could choose $T_i$ so that $Y^{[h]}$ is the reduced closed subscheme of $X$ associated to the pull-back of the closed subscheme of $\Spec A_0$ defined by $(\fkp,T_1,\ldots, T_{n-1-h})$, i.e.\ the closure of the set of closed points where the reduction of $\Sigma^\univ$ has height $\geq n-h$).

Here we explain a proof using the integral model of unitary Shimura varieties (we have already recalled in Section 2.4 that this gives the existence of $U$ and $x$). This has an advantage that it spells out the corresponding blow-ups of Shimura varieties, and calculations in the last section imply that we have a generalized semistable model of Shimura varieties in this case.

We use the notation from \cite{HT}, especially Chapter III.4. The integral model $X_{U^p,m}$ of Shimura varieties defined in p.109 of \cite{HT} is a proper flat scheme of relative dimension $n-1$ over $\mO_{F,w}=\mO$, and its special fiber $\ol{X}_{U^p,m}=X_{U^p,m}\times_{\Spec \mO}\Spec k$ has a stratification by reduced closed subschemes $\ol{X}_{U^p,m}^{[h]}$ of dimension $h$ for $0\leq h\leq n-1$. This $\ol{X}_{U^p,m}^{[h]}$ is the closure of the set of closed points where the associated 1-dimensional Barsotti-Tate $\mO$-module $\mG$ has etale height $\leq h$ (defined in p.111 of \cite{HT}). 

We fix a multi-index $m$ with $m_1=1$, and denote $X_{U^p,m}$ simply by $\Sh$, and similarly $\ol{\Sh}:=\ol{X}_{U^p,m},\ \ol{\Sh}^{[h]}:=\ol{X}_{U^p,m}^{[h]}$ etc. We recalled in Section 2.4 that the completion $\widehat{\mO}_{\Sh,\ol{s}}$ of the strict local ring of $\Sh$ at any geometric closed point $\ol{s}$ centered in $s\in \ol{\Sh}^{[0]}$ is isomorphic to our deformation ring $A$, i.e.\ $X=\Spec A\cong \Spec \widehat{\mO}_{\Sh,\ol{s}}$. We fix the resulting morphism
\[X\cong \Spec \widehat{\mO}_{\Sh,\ol{s}}\lra \Sh.\]

\blem \label{yhshh}
We have $Y^{[h]}=X\times_\Sh\ol{\Sh}^{[h]}$.
\elem

\bprf
We can show directly that $Y^{[h]}$ is the reduced subscheme of $X$ associated to the locus where the deformation of $\Sigma_n$ has height $\geq n-h$. Here we give a proof closer in spirit to how we defined $Y^{[h]}$. In \cite{Man} 3.2 (see Remark 10(2)), the following decomposition of $\ol{\Sh}^{[h]}$ is constructed (their intersections with the open strata are defined in \cite{HT}, p.115, denoted $\ol{X}_{U^p,m,M}$). If $M$ is a free $k=\mO/\fkp$-submodule of $\fkp^{-1}\Lambda_{11}/\Lambda_{11}$ of rank $n-h$ (here $\Lambda_{11}\cong (\mO^n)^\vee$ is a free $\mO$ module of rank $n$ used in the definition of the Shimura variety), then $\ol{\Sh}_M^{[h]}$ is the unique closed subscheme of $\Sh$ such that, for any scheme $T$ over $\Sh$, a $T$-valued point $T\ra \Sh$ factors through $\ol{\Sh}_M^{[h]}$ if and only if the base change of the Drinfeld level $\fkp$ structure $\alpha$ on the universal Barsotti-Tate $\mO$-module $\mG$ on $\Sh$ of height $n$:
\[\alpha_T:(\fkp^{-1}\Lambda_{11}/\Lambda_{11})_T\lra \mG[\fkp]_T\]
vanishes on $M_T$. Then $\ol{\Sh}^{[h]}=\bigcup_M\ol{\Sh}_M^{[h]}$, where $M$ runs through all $k$-subspaces of $\fkp^{-1}\Lambda_{11}/\Lambda_{11}$ of rank $n-h$.

Now, the universal formal $\mO$-module $\Sigma^\univ\otimes A$ and the universal level $\fkp$ structure $\phi$ are the pull-backs under our morphism $X\ra \Sh$ of the Barsotti-Tate $\mO$-module $\mG/\Sh$ and the level $\fkp$ structure $\alpha$. Therefore, for any scheme $T\ra \Sh$ pulled back to $T_X\ra X$, the condition for $T_X$ to factor through $X\times_\Sh\ol{\Sh}_M^{[h]}$ translates as follows. Expressing the elements of $(\fkp^{-1}\Lambda_{11}/\Lambda_{11})_{T_X}\cong (k^n)^\vee$ by the standard basis, vanishing of $\alpha_{T_X}$ on $M_{T_X}$ is written as:
\[(a_1,\ldots,a_n)\in M\Lra [a_1](X_1)+_{\Sigma^\univ}\cdots +_{\Sigma^\univ}[a_n](X_n)=0.\]
This RHS is exactly the defining equation of $P_{\ul{a}}$. If we take a basis $\ul{a}_1,\ldots,\ul{a}_{n-h}$ of $M$ and define the hyperplanes $M_i:=M_{\ul{a}_i}$ of $\Ps$, then $N:=\bigcap_{i=1}^{n-h}M_i\subset \Ps$ is a $k$-rational linear subspace of $\Ps$ of dimension $h-1$. We see that $X\times_\Sh\ol{\Sh}_M^{[h]}$ is the closed subscheme of $X$ defined by the ideal $(P_{\ul{a}_1},\ldots,P_{\ul{a}_{n-h}})$, namely $Y_N$ by Lemma \ref{yn}. As $M$ runs through the rank $n-h$ free submodules of $(k^n)^\vee$, clearly $N$ runs through all $k$-rational linear subspaces of $\Ps$ of dimension $h-1$.
\eprf

\bde
We denote the blow-up of $\Sh$ at $\ol{\Sh}^{[0]}$ by $\Sh_1$, and inductively define $\Sh_{h+1}\ra \Sh_h$ as the blow-up of $\Sh_h$ at the proper transform of $\ol{\Sh}^{[h]}$ inside $\Sh_i$, for $1\leq h\leq n-2$.
\ede

\blem
We have canonically $Z_\st=Z_{n-1}\cong X\times_\Sh\Sh_{n-1}$.
\elem

\bprf
We show that $Z_h\ra \Sh_h$ is flat and $Z_h\cong X\times_\Sh\Sh_h$ by induction. Corollary \ref{blcp} shows that it is true for $h=1$, the strict henselisation and the completion being flat. The above lemma shows that the center of blow-up for $Z_{h+1}\ra Z_h$ is the pull back of the center of blow-up for $\Sh_{h+1}\ra \Sh_h$. Therefore we have $Z_{h+1}\cong Z_h\times_{\Sh_h}\Sh_{h+1}$ by Lemma\ \ref{blbc}, therefore $Z_{h+1}\cong X\times_\Sh\Sh_{h+1}$ by our inductive hypothesis.
\eprf

Now if we denote the base change of $\Sh,\Sh_{n-1}$ by $S=\Spec W\ra \Spec \mO$ by $\Sh_S,\Sh_\st$, we have a diagram of $S$-schemes
\[\xymatrix{
Y_{\Ps,\st} \ar[r] \ar[d]_\cong & Z_\st \ar[r] \ar[d]_f & X \ar[d] \ar[dr]\\
Y_{\ol{s},\st} \ar[r] & \Sh_\st \ar[r] & \Sh_S \ar[r] & S
}\]
whose second square is cartesian by the above lemma. We can regard the geometric closed point $\ol{s}\ra s\in \Sh$ that we chose before as the closed point of $\Sh_S$, and the inverse image of the closed point $x\in X$ under $Z_\st\ra X$, namely the unique proper component $Y_{\Ps,\st}$ of $Z_\st\times_{\Spec W}\Spec \ol{k}$, is isomorphic under $f:Z_\st\ra \Sh_\st$ to the inverse image $Y_{\ol{s},\st}$ of $\ol{s}$ under $\Sh_\st\ra \Sh_S$. 

\blem \label{compisom}
For all closed points $z\in Y_{\Ps,\st}\subset Z_\st$, the local $W$-homomorphism between the complete local rings $\widehat{\mO}_{\Sh_\st,f(z)}\ra \widehat{\mO}_{Z_\st,z}$ induced by $f$ is an isomorphism.
\elem

\bprf
The morphism $\Sh_\st\ra \Sh_S$ maps $f(z)$ to $\ol{s}$, hence it induces a local homomorphism between complete local rings $\widehat{\mO}_{\Sh_S,\ol{s}}\ra \widehat{\mO}_{\Sh_\st,f(z)}$, or $\Spec \widehat{\mO}_{\Sh_\st,f(z)}\ra X$. This gives $\Spec \widehat{\mO}_{\Sh_\st,f(z)}\ra Z_\st\cong X\times_{\Sh_S}\Sh_\st$ whose image of the closed point must be $z$ as $f^{-1}(f(z))=\{z\}$, therefore the desired inverse $\widehat{\mO}_{Z_\st,z}\ra \widehat{\mO}_{\Sh_\st,f(z)}$ of the local homomorphism in the proposition.
\eprf

Therefore we have:

\bpr \label{shgsst}
\benu
\item The proper flat $S$-scheme $\Sh_\st$ has generalized semistable reduction at the points in $Y_{\ol{s},\st}$ of the special fiber, in the sense that for any $x\in Y_{\ol{s},\st}$ there exists an etale morphism $U\ra \Sh_\st$ with $x\in \im U$ and an etale $S$-morphism:
\[\phi_U:U\lra \Spec W[T_1,\ldots,T_n]/(T_1^{e_1}\cdots T_d^{e_d}-\varpi)\ \ (d\leq n),\]
where integers $e_i$ are all prime to $\chara k$.
\item The canonical base change morphisms of nearby cycle sheaves $f^*(R^i\psi\Lambda|_{Y_{\ol{s},\st}})\lra R^i\psi\Lambda|_{Y_{\Ps,\st}}$ (for $\Lambda=\ol{\Q}_\ell$) are isomorphisms for all $i$.
\eenu
\epr

\bprf
(i) This follows from Lemma \ref{compisom}, using the fact that we can characterize the generalized semistable reduction by looking at the completion of the strict local rings. From the lack of immediate reference, we give a sketch of its proof. Assume that the complete local ring at a closed point $x$ of the scheme $X$ of finite type over $S$ is isomorphic to $W[[T_1,\ldots,T_n]]/(T_1^{e_1}\cdots T_d^{e_d}-\varpi)\ (d\leq n)$, where the integers $e_i$ are all prime to $\chara k$. Then it is a regular scheme with its reduced special fiber being a normal crossing divisor with all the multiplicities prime to $\chara k$ (i.e.\ the generators $t_1,\ldots,t_d$ of the components passing through $x$ are the part of a regular system of parameters). Therefore, some etale neighborhood of $x$ in $X$ is a regular scheme with its reduced special fiber being a normal crossing divisor with all the multiplicities prime to $\chara k$. Then by sending $T_i$ to the local generators $t_i$ of the divisors crossing at $x$ we obtain an unramified map from an etale neighborhood $\Spec A$ of $x$ in $X$ to $\Spec W[T_1,\ldots,T_n]$. By \cite{EGA} IV, (18.4.7), we can decompose this morphism into an etale morphism $\Spec B\ra \Spec W[T_1,\ldots,T_n]$ and a closed immersion $\Spec A\ra \Spec B$. Hence $A=B/I$ for an ideal $I$ of $B$, and the inverse image of $I$ in $W[T_1,\ldots,T_n]$ must be of the form $(T_1^{e_1}\cdots T_d^{e_d}-u\varpi)$ with a unit $u$. This induces an etale map $\Spec A\ra \Spec W[T_1,\ldots,T_n]/(T_1^{e_1}\cdots T_d^{e_d}-u\varpi)$, and as all the $e_i$'s are prime to $p$, we can take $u=1$ etale locally, hence we obtain the desired etale map.

(ii) This follows from the regular base change theorem for nearby cycle sheaves (\cite{Fu}, Cor.\ 7.1.6), as $f$ is regular (because $\Sh_S$ is excellent and $X\ra \Sh_S$ is regular). (Or, using Lemma \ref{compisom}, the morphisms induced on the stalks are isomorphisms by the formal invariance theorem (\cite{Fu}, Cor.\ 7.1.7), i.e.\ the isomorphism $(R^i\psi\Lambda)_z\cong H^i(\Spec(\widehat{\mO}_z)_{\ol{\eta}},\Lambda)$, where $_{\ol{\eta}}$ denotes the geometric generic fiber.)
\eprf

\brem
{\rm We can show that the whole $\Sh_\st$ has generalized semistable reduction in a similar way.}
\erem

\section{A model $Z_n$ containing the Deligne-Lusztig variety}

\subsection{Base change and normalization}

Now we proceed to construct a model of $X$ over a tamely ramified extension of $W$, namely $W_n:=W(\varpi_n)$ where $\varpi_n:=\varpi^{1/(q^n-1)}$, which includes a Deligne-Lusztig variety inside the special fiber. This is done simply by taking the base change of $Z_1$ over $S_n:=\Spec W_n$ and normalizing.

\bde
\benu
\item Let $Z_n$ be the normalization of $Z_1\times_SS_n$.
\item Let $U_n$ and $Y_n$ respectively be the inverse images of $Y_\Ps^{(n)}$ and $Y_\Ps$ under the morphism $Z_n\ra Z_1$. The $\ol{k}$-scheme $Y_n$ is a proper subvariety of the special fiber of $Z_n$, and $U_n$ is an open subvariety of $Y_n$.
\eenu
\ede

First we define an open affine subscheme $\Spec C_1$ of $\Spec B_1$ (the completion of an affine open of $Z_1$ along the exceptional divisor), which has $Y_\Ps^{(n)}$ as the special fiber. As $Y_\Ps^{(n)}$ is the complement in $Y_\Ps$ of the intersection with all the other irreducible components, namely $Y_{M,1}$ for all $M$, we only need to invert the equations which reduces $\bmod \fkp$ to the defining equations of $Y_{M,1}$. We start from the affine subscheme defined in the proof of Proposition \ref{z0eq}, namely:
\[\Spec B_1=W[V_1,\ldots,V_{n-1}][[\widetilde{X}_n]] \bigm/ \bigl(\tilde{u}\cdot \widetilde{X}_n^{q^n-1}\cdot \!\!\!\prod_{\ul{a}\in k^n\setminus\{\ul{0}\}}P'_{\ul{a}}-\varpi\bigr).\]
with $\tilde{u}\in 1+(\widetilde{X}_n)\subset W[V_1,\ldots,V_{n-1}][[\widetilde{X}_n]]^\times$. The exceptional divisor inside $\Spec (B_1\otimes_W\ol{k})$ is the affine space $\A^{n-1}$, the complement of the hyperplane $X^*_n=0$ in $\Ps$. Now in order to remove all the $k$-rational hyperplanes $\sum_{i=0}^{n-1}\tilde{a_i}V_i+a_n=0$ of $\A^{n-1}$ from the special fiber, we invert the element $\prod_{\ul{a}\in k^n\setminus\{\ul{0}\}}P'_{\ul{a}}$ and define:
\[C_1:=B_1\left[\frac{1}{\prod P'_{\ul{a}}}\right]=W[V_1,\ldots,V_{n-1}][[\widetilde{X}_n]]\left[\frac{1}{\prod P'_{\ul{a}}}\right]\bigm/\bigl(\tilde{u}\cdot \widetilde{X}_n^{q^n-1}\cdot \prod P'_{\ul{a}}-\varpi\bigr),\]
because $Y_{M,1}$ for $M=M_{\ul{a}}$ was defined by $P'_{\ul{a}}=0$.

Now we describe the equation of the variety $U_n$, which is achieved by normalizing the ring $C_1\otimes_WW_n$, which is:
\[W_n[V_1,\ldots,V_{n-1}][[\widetilde{X}_n]]\left[\frac{1}{\prod P'_{\ul{a}}}\right]\bigm/(\tilde{u}\cdot \widetilde{X}_n^{q^n-1}\cdot\prod P'_{\ul{a}}-\varpi_n^{q^n-1}).\]
For this we adjoin the element $V_n:=\varpi_n/\widetilde{X}_n$ of its total quotient ring to this ring, as the equation:
\[\varpi_n^{q^n-1} = \tilde{u}\cdot \widetilde{X}_n^{q^n-1}\cdot \prod P'_{\ul{a}},\ \ \text{i.e.}\ \ (\varpi_n/\widetilde{X}_n)^{q^n-1}=\tilde{u}\cdot \prod P'_{\ul{a}}\]
shows that $V_n=\varpi_n/\widetilde{X}_n$ is integral over $C_1\otimes_WW_n$, and moreover it is a unit, as $V_n^{q^n-1}=\tilde{u}\cdot \prod P'_{\ul{a}}$ is inverted. By using the fact that inverting $\prod P'_{\ul{a}}$ is equivalent to inverting $V_n$, we compute the ring $C_1\otimes_WW_n[V_n]$ as follows:
\begin{align*}
 & W_n[V_1,\ldots,V_n]\left[\frac{1}{V_n}\right][[\widetilde{X}_n]]\bigm/(\varpi_n-V_n\widetilde{X}_n,\ \tilde{u}\cdot \widetilde{X}_n^{q^n-1}\prod P'_{\ul{a}}-\varpi_n^{q^n-1})\\
 &= W_n[V_1,\ldots,V_n]\left[\frac{1}{V_n}\right][[\widetilde{X}_n]]\bigm/(\varpi_n-V_n\widetilde{X}_n,\ \widetilde{X}_n^{q^n-1}(\tilde{u}\cdot \prod P'_{\ul{a}}-V_n^{q^n-1})).\\
\end{align*}
Now we claim the following:

\bpr
Consider the spectrum $\Spec C_n$ of the ring:
\[C_n=W_n[V_1,\ldots,V_n]\left[\frac{1}{V_n}\right][[\widetilde{X}_n]]\bigm/(\varpi_n-V_n\widetilde{X}_n,\ \tilde{u}\cdot \prod P'_{\ul{a}}-V_n^{q^n-1}),\]
which is a closed subscheme of $\Spec (C_1\otimes_WW_n[V_n])$.
\benu
\item $\Spec C_n$ is regular.
\item $\Spec C_n$ is the normalization of $C_1\otimes_WW_n$. 
\eenu
\epr

\bprf
(i) It can be seen by applying Jacobian computation directly, but here we prove it by showing that $\Spec C_n$ is formally smooth over $S_n$. The generic fiber is apparently formally smooth as it is unchanged from that of $\Spec B_1\times_SS_n$ and $\Spec B_1$ was the completion of an affine open of $Z_1$, which in turn had the same generic fiber as $X$. The special fiber of $\Spec C_n$ (equal to the locus of $\widetilde{X}_n=0$, as $V_n=\varpi_n/\widetilde{X}_n$ is a unit), i.e.\ the spectrum of the ring:
\[C_n\otimes_{W_n}\ol{k}=\ol{k}[V_1,\ldots,V_n]\left[\frac{1}{V_n}\right]\bigm/ \Big(\prod (P'_{\ul{a}}\bmod \widetilde{X}_n) -V_n^{q^n-1}\Big)\]
(here we used $\tilde{u}\in 1+(\widetilde{X}_n)$), is a smooth affine variety, because we have
\[\prod (P'_{\ul{a}}\bmod \widetilde{X}_n) = \prod_{\ul{a}\in k^n\setminus\{\ul{0}\}}(a_1V_1 + \cdots + a_{n-1}V_{n-1}+a_n)\]
by (\ref{z0eq1}). Now it remains to show that $\Spec C_n$ is flat over $S_n$. Starting from the flat $W_n$-algebra $W_n[V_1,\ldots,V_n][1/V_n][[\widetilde{X}_n]]$, we will apply the following lemma twice:

\blem
Let $(A,\fkm)$ be a noetherian local ring, and let $B$ be a noetherian flat $A$-algebra such that every maximal ideal of $B$ lies over $\fkm$. If $b\in B$ is $B/\fkm B$-regular (i.e.\ multiplication by $b$ on $B/\fkm B$ is injective), then $B/(b)$ is flat over $A$.
\elem

\bprf
It is a special case of \cite{Mat}, Th.\ 22.6.
\eprf

Letting $A=W_n$, first for $B=W_n[V_1,\ldots,V_n][1/V_n][[\widetilde{X}_n]]$, we see that $B/\varpi_nB=\ol{k}[V_1,\ldots,V_n][1/V_n][[\widetilde{X}_n]]$ is an integral domain, where $\varpi_n-V_n\widetilde{X}_n\in B$ does not reduce to zero, hence we have the flatness of $B/(\varpi_n-V_n\widetilde{X}_n)$. 

Next, letting $B=W_n[V_1,\ldots,V_n][1/V_n][[\widetilde{X}_n]]\bigm/(\varpi_n-V_n\widetilde{X}_n)$, and seeing that 
\[B/\varpi_nB=\ol{k}[V_1,\ldots,V_n]\left[\frac{1}{V_n}\right][[\widetilde{X}_n]]\bigm/(V_n\widetilde{X}_n)=\ol{k}[V_1,\ldots,V_n]\left[\frac{1}{V_n}\right]\]
is again an integral domain where $\tilde{u}\cdot \prod P'_{\ul{a}}-V_n^{q^n-1}\in B$ does not reduce to zero, flatness of $C_n$ over $W_n$ follows as desired.

(ii) Clearly the generic point of $C_n$ coincides with that of $C_1\otimes_WW_n[V_n]$, i.e.\ that of $C_1\otimes_WW_n$, and $C_n$ is finite over $C_1\otimes_WW_n$. Hence the assertion follows from (i).
\eprf

Now $\Spec (C_1\otimes_WW_n)$ is an affine open of $\Spec B_1\times_SS_n$, which is the completion along the exceptional divisor of an affine open of $Z_1\times_SS_n$. Therefore $\Spec C_n$ is an affine open of the completion along the exceptional divisor of the normalization $Z_n$. Moreover, as the special fiber of $\Spec C_1$ is naturally isomorphic to $Y_\Ps^{(n)}$, which is contained in the exceptional divisor, the special fiber $\Spec (C_n\otimes_{W_n}\ol{k})$ must be isomorphic to the inverse image of $Y_\Ps^{(n)}$ inside $Z_n$, namely $U_n$.

By changing the coordinates on $U_n$ as follows:
\[X'_i:=V_i/V_n\ \ (1\leq i\leq n-1),\ \ \ X'_n:=1/V_n\]
we have the following description of $U_n\cong \Spec (C_n\otimes_{W_n}\ol{k})$:
\begin{align*}
U_n &= \Spec \ol{k}[V_1,\ldots,V_n]\left[\frac{1}{V_n}\right]\Bigm/ \Big(\prod_{\ul{a}\in k^n\setminus\{\ul{0}\}}(a_1V_1+\cdots +a_{n-1}V_{n-1}+a_n)-V_n^{q^n-1}\Big)\\
&= \Spec \ol{k}[X'_1,\ldots,X'_n] \left[\frac{1}{X'_n}\right] \Bigm/ \Big(\prod_{\ul{a}\in k^n\setminus\{\ul{0}\}}\bigg(a_1\frac{X'_1}{X'_n}+\cdots +a_{n-1}\frac{X'_{n-1}}{X'_n}+a_n\bigg) - \bigg( \frac{1}{X'_n} \bigg)^{q^n-1} \Big)\\
&= \Spec \ol{k}[X'_1,\ldots,X'_n]\left[\frac{1}{X'_n}\right]\Bigm/ \Big(\prod_{\ul{a}\in k^n\setminus\{\ul{0}\}}(a_1X'_1+\cdots +a_nX'_n)-1\Big)
\end{align*}

Note that, by unwinding the definitions, we see that the coordinates $X'_1,\ldots,X'_n$ are related to the original $\widetilde{X}_1,\ldots,\widetilde{X}_n$ through $X'_i=\widetilde{X}_i/\varpi_n$ for each $1\leq i\leq n$.

\bpr \label{equ2}
The equation of $U_n$ is given by:
\[U_n=\Spec \ol{k}[X'_1,\ldots,X'_n]\Bigm/ \Big(\prod_{\ul{a}\in k^n\setminus\{\ul{0}\}}(a_1X'_1+\cdots +a_nX'_n)-1\Big)\]
which is a smooth affine variety over $\ol{k}$.
\epr

\subsection{Action of $GL_n(k)$ and the inertia group}

Summing up, we have following diagram of schemes:
\[
\xymatrix{
U_n \ar[r]\ar[d] & Y_n \ar[r]\ar[d] & Z_n \ar[r]\ar[d] & S_n \ar[d] \\
Y_\Ps^{(n)} \ar[r] & Y_\Ps \ar[r] & Z_1 \ar[r] & S
}\]
Here all the vertical maps are surjective. The $W$-schemes $Z_1$ and $Z_n$ are normal models of $X$ and $X\times_SS_n$ respectively. The $\ol{k}$-schemes $Y_\Ps$ and $Y_n$ are the unique proper components (with multiplicity $q^n-1,1$ respectively) of their special fibers, and $Y_\Ps^{(n)}$ and $U_n$ are open affine subvarieties of them. Also:
\benu
\item The right action of $GL_n(k)$ on $X$ extend to all the schemes in the above diagrams, and all maps are $GL_n(k)$-equivariant.
\item The special fiber of $Z_n$, and therefore also $Y_n$ and $U_n$, has the geometric inertia action of $I_K$, which factors through the finite quotient:
\[ I_K\ra \gal(\Frac W_n/\Frac W)\cong k_n^\times. \]
\eenu
Now we describe these actions on $U_n$ explicitly by the coordinates introduced above. We denote the image of $X'_i$ in the coordinate ring of $U_n$ by $\ol{X'_i}$.

\bpr \label{actionun}
\benu
\item The group $GL_n(k)$ acts on $U_n$ via right linear action on the row vector $(\ol{X'_1},\ldots,\ol{X'_n})$.
\item The covering $f:U_n\ra Y_\Ps^{(n)}$ induced from the finite map $Z_n\ra Z_1$ is a finite Galois etale covering with the automorphism group $\gal(\Frac W_n/\Frac W)$, which is canonically isomorphic to $\MU_{q^n-1}\cong k_n^\times$. Its action on the coordinate ring is described by $\ol{X'_i}\mapsto \zeta^{-1} \ol{X'_i}$ for all $1\leq i\leq n$, where $\zeta\in \MU_{q^n-1}$ denotes the image under the above canonical isomorphism. \eenu
\epr

\bprf
(i) This is easily seen by substituting $\widetilde{X}_j=\varpi_nX'_j$ to the original action (\ref{glaction}) and using Proposition \ref{paiprop}(i) to see that the action of $(a_{ij})\in GL_n(k) \bmod \varpi_n^2$ is
\[\varpi_n\ol{X'_j}\longmapsto a_{1j}(\varpi_n\ol{X'_1})+\cdots +a_{nj}(\varpi_n\ol{X'_n})\]
and dividing both sides by $\varpi_n$.

(ii) This is also clear by the relation $X'_i=\widetilde{X}_i/\varpi_n$ for each $1\leq i\leq n$, and the inertia action on $S_n$ being defined by $\varpi_n\mapsto \zeta\varpi_n$ for $\zeta\in \MU_{q^n-1}\cong \gal(\Frac W_n/\Frac W)$.
\eprf


We will observe in the next chapter that the finite etale covering $U_n\ra Y_\Ps^{(n)}$ of smooth affine varieties over $\ol{k}$, with the right action of $GL_n(k)\times I_K$, is isomorphic to the Deligne-Lusztig variety for $GL_n(k)$ and its maximally non-split torus $T$ with $T(k)\cong k_n^\times$, defined in \cite{DL}, Sections 2.1 and 2.2.

\section{Vanishing cycles}

\subsection{Review on vanishing cycle cohomology}

Here we recall the method of calculating vanishing cycle cohomology of strict local rings over strict henselian discrete valuation ring via the special fiber of a model. A similar argument can be found in \cite{Br}.

In this subsection, let $K$ be a strict henselian discrete valuation field, let $\mO$ be its ring of integers and $k$ be its residue field with $\chara k=p>0$. We denote the generic point and the closed point of $S:=\Spec \mO$ by $\eta:=\Spec K$ and $s:=\Spec k$, and the geometric generic point by $\ol{\eta}=\Spec \ol{K}$, where $\ol{K}$ is a separable closure of $K$. We let $\ol{S}:=\Spec \mO_{\ol{K}}$, where $\mO_{\ol{K}}$ is the integral closure of $\mO$ in $\ol{K}$, which is a non-discrete valuation ring. The generic point and the closed point of $\ol{S}$ are respectively $\ol{\eta}$ and $s$.

We consider an $n$-dimensional strict henselian local $\mO$-algebra $A$ with local structure homomorphism $\mO\ra A$, and let $X:=\Spec A$. We denote the base change of $X/S$ over $s,\eta,\ol{\eta},\ol{S}$ respectively by $X_s,X_\eta,X_{\ol{\eta}},\ol{X}$. The above schemes fit into the following diagrams:
\[\xymatrix{
s\ar@{=}[d] \ar[r] & \ol{S} \ar[d] & \ol{\eta} \ar[d] \ar[l]\\
s\ar[r]  & S & \eta \ar[l]
}\ \ \ \ \ 
\xymatrix{
X_s\ar@{=}[d] \ar[r]^-{\ol{i}} & \ol{X} \ar[d] & X_{\ol{\eta}} \ar[d] \ar[l]_-{\ol{j}}\\
X_s\ar[r]^-i & X & X_\eta \ar[l]_-j
}\]
We are interested in the $\ell$-adic etale cohomology groups $H^*(X_{\ol{\eta}},\ol{\Q}_\ell)$ of the $(n-1)$-dimensional affine scheme $X_{\ol{\eta}}$ over $\ol{\eta}=\Spec \ol{K}$, where $\ell$ is a prime not equal to $p$. Hereafter we denote the constant $\ell$-adic sheaf $\ol{\Q}_\ell$ simply by $\Lambda$. By the Leray spectral sequence for $\ol{j}$, we have canonical isomorphisms:
\[H^i(X_{\ol{\eta}},\Lambda)\cong \bH^i(\ol{X},R\ol{j}_*\Lambda)\cong (R^i\ol{j}_*\Lambda)_x,\]
where $\bH$ denotes the hypercohomology and $R\ol{j}_*\Lambda$ is the object in the derived category of $\ell$-adic sheaves on $\ol{X}$, and $x$ denotes the unique closed point of $X$ and $\ol{X}$. We denoted by $(R^i\ol{j}_*\Lambda)_x$ the stalk at $x$ of the $\ell$-adic sheaf $R^i\ol{j}_*\Lambda$ over $\ol{X}$. Constructibility of the $\ell$-adic sheaf $R\ol{j}_*\Lambda$ (from which follows that $H^i(X_{\ol{\eta}},\Lambda)$ has finite rank) follows from \cite{SGA4.1/2}, [Finitude], under the hypothesis that $X$ is the strict local ring of a scheme locally of finite type over $S$.

We can try to compute these cohomology groups by constructing good models of $X$, obtained by repeated blow-ups and normalizations over ramified extensions of $K$. 

Let $K'$ be a finite extension of $K$ and $S':=\Spec \mO'$ be the normalization of $S$ in $K'$. Let $X':=X\times_SS',\ \eta':=\Spec K',\ X'_{\eta'}:=X'\times_{S'}\eta',\ j':X'_{\eta'}\ra X'$ and let $x$ denote also the closed point of $X'$. Suppose we have a proper morphism $p:Z\ra X'$ over $S'$ which is an isomorphism on the generic fibers (i.e.\ $p|_{Z_{\eta'}}=\id$ where $Z_{\eta'}:=Z\times_{S'}\eta'$), and denote the inverse image of $x$ by $Y:=p^{-1}(x)$. Let $i_Z,j_Z$ be the inclusions $Y\ra Z$ and $X'_{\eta'}\ra Z$ respectively, and we denote the objects obtained from $p,Z,i_Z,j_Z$ by the base change under $\ol{S}\ra S'$ respectively by $\ol{p}, \ol{Z}, \ol{i_Z}, \ol{j_Z}$. Note that $Y=\ol{p}^{-1}(x)$ is a proper subscheme of the special fiber $Z_s:=Z\times_Ss$ of $Z$. We can describe the situation by the following diagrams:
\[\xymatrix{
Y\ar[d] \ar[r]^-{i_Z} & Z \ar[d]^-p & X'_{\eta'} \ar@{=}[d] \ar[l]_-{j_Z}\\
x\ar[r]  & X' & X'_{\eta'} \ar[l]_-{j'}
}\ \ \ \ \ 
\xymatrix{
Y\ar[d] \ar[r]^-{\ol{i_Z}} & \ol{Z} \ar[d]^-{\ol{p}} & X_{\ol{\eta}} \ar@{=}[d] \ar[l]_-{\ol{j_Z}}\\
x\ar[r]  & \ol{X} & X_{\ol{\eta}} \ar[l]_-{\ol{j}}
}\]
Then, if we denote by $R\psi\Lambda$ the {\it nearby cycle sheaves} $\ol{i_Z}^*R\ol{j_Z}_*\Lambda$ in the derived category of $\ell$-adic sheaves on $Y$, the proper base change theorem (\cite{SGA4} Expos\'e XIII) gives the following canonical isomorphisms:
\[(R\ol{j}_*\Lambda)_x\cong (R\ol{p}_*R\ol{j_Z}_*\Lambda)_x\cong R\Gamma(Y,R\psi\Lambda).\]
Therefore we have:

\bpr \label{specseq}
We have a canonical isomorphism $H^i(X_{\ol{\eta}},\Lambda)\cong \bH^i(Y,R\psi\Lambda)$. (As it is canonical, if a group $G$ acts on $X$ and $Z\ra X$ is $G$-equivariant, then this isomorphism is $G\times I_K$-equivariant. This holds for similar isomorphisms in what follows.)
\epr

We will also need the following consequence of the proper base change theorem:

\bpr \label{finvancyc}
Let $K'_1\subset K'_2$ be successive finite extensions of $K$, and suppose we have the above situation over each field: 
\[\xymatrix{
Y_1\ar[d] \ar[r]^-{i_{Z_1}} & Z_1 \ar[d]^-{p_1} & X'_{\eta'_1} \ar@{=}[d] \ar[l]_-{j_{Z_1}}\\
x\ar[r]  & X'_1 & X'_{\eta'_1} \ar[l]_-{j'_1}
}\ \ \ \ \ 
\xymatrix{
Y_2\ar[d] \ar[r]^-{i_{Z_2}} & Z_2 \ar[d]^-{p_2} & X'_{\eta'_2} \ar@{=}[d] \ar[l]_-{j_{Z_2}}\\
x\ar[r]  & X'_2 & X'_{\eta'_2} \ar[l]_-{j'_2}
}\]
with a proper morphism $f:Z_2\ra Z_1$ which induces the isomorphism $X'_{\eta'_1}\times \eta'_2\cong X'_{\eta'_2}$. If we denote the nearby cycle sheaves on $Y_1,Y_2$ respectively by $R\psi_1\Lambda, R\psi_2\Lambda$, we have a canonical isomorphism $R\psi_1\Lambda\cong Rf_*R\psi_2\Lambda$. In particular, if $f$ is finite, we have $R^i\psi_1\Lambda\cong f_*R^i\psi_2\Lambda$ for each $i$.
\epr

Lastly, in order to apply the results from Section 4 to compute the vanishing cycle cohomology, we need the results of T.\ Saito on the nearby cycle sheaves on generalized semistable schemes, following \cite{Sa} (see also \cite{RZ}).

Let $X$ be a scheme locally of finite type over the strict henselian trait $S=\Spec \mO$. We denote the inclusions $X_s\ra X$ and $X_\eta\ra X$ respectively by $i$ and $j$, and denote by $\ol{X},X_{\ol{\eta}},\ol{i},\ol{j}$ the objects obtained by the base change under $\ol{S}\ra S$ of the corresponding objects. Consider the nearby cycle sheaves $R^q\psi\Lambda:=\ol{i}^*R^q\ol{j}_*\Lambda$, which are constructible $\ell$-adic sheaves on $X_s$ with the action of the inertia group $I_K=\gal(\ol{K}/K)$ of $K$.

Assume that $X$ has generalized semistable reduction in the sense of Proposition \ref{shgsst}(i). Let $\{Y_i\}_{i\in I}$ be the irreducible components of the special fiber $X_s$, and for any finite subset $J\subset I$, we denote $Y_J:=\bigcap_{i\in J}Y_i$ and $Y_J^0:=Y_J\setminus \bigcup_{i\notin J}Y_i$. Let $e_i$ denote the multiplicity of $Y_i$ in $X_s$ for $i \in I$. By Proposition 6 of \cite{Sa}, we have:

\bpr \label{gensstvan}
Let $d=d_J$ be the greatest common divisor of $\{e_i\}_{i\in J}$, and $K_d$ be the unique tamely ramified extension over $K$ of degree $d$. Then we have the following canonical isomorphisms as constructible $\ell$-adic sheaves with $I_K$-action:
\begin{align*}
R^0\psi\Lambda|_{Y_J^0} &\cong \ind_{I_{K_d}}^{I_K}\Lambda\ \ (\text{etale locally}),\\
R^q\psi\Lambda|_{Y_J^0} &\cong R^0\psi\Lambda|_{Y_J^0}\otimes_{\Z_\ell} \bigwedge^q \Hom \Big(\Ker\bigl(\bigoplus_{i\in J}\Z_\ell\ra \Z_\ell\bigr)(1),\ \Lambda \Big),
\end{align*}
where the map $\bigoplus_{i\in J}\Z_\ell\ra \Z_\ell$ is defined by $1_i\mapsto e_i$, and $(1)$ denotes the Tate twist.
\epr

\bcor \label{corsaito}
\benu
\item Each $R^q\psi\Lambda$ is locally constant on $Y_J^0$ of rank equal to $d\cdot \binom{|J|-1}{q}$, where $\binom{|J|-1}{q}=0$ for $q>|J|-1$ by convention. Moreover, the inertia group $I_K$ acts on any $R^q\psi\Lambda$ through the finite cyclic Galois group $I_K/I_{K_d}\cong \MU_d$. 
\item {\rm (Cor.\ 1 to Prop.\ 6 of \cite{Sa})} On $Y_J^0$ with $|J|\neq 1$, in the Grothendieck group of smooth $\ell$-adic sheaves with $I_K$-action, the following alternating sum vanishes:
\[[R\psi\Lambda|_{Y_J^0}]=\sum_q(-1)^q[R^q\psi\Lambda|_{Y_J^0}]=0.\]
\eenu
\ecor


\subsection{Review of the Deligne-Lusztig theory}

Here we recall some results of the Deligne-Lusztig theory \cite{DL} that are relevant to our case (see also \cite{Se}).

Let $G$ be a connected reductive group defined over a finite field $k$, and $F:G\ra G$ be the Frobenius morphism. For an $F$-stable maximal torus $T$ and a Borel subgroup $B$ containing $T$, Deligne and Lusztig construct a $G^F$-equivariant finite etale Galois covering of smooth varieties over $\ol{k}$ with right $G^F$-actions:
\[f:\widetilde{X}_{T\subset B}\ra X_{T\subset B},\]
with Galois group $T^F$. By decomposing the $\ell$-adic sheaf as $f_*\ol{\Q}_\ell\cong \bigoplus_\theta\mF_\theta$, where $\theta$ runs through the characters of $T^F$, we define a virtual $G_F$-representation
\[R^\theta_T:=\sum_i(-1)^iH^i_c(X_{T\subset B},\mF_\theta),\]
which turns out to be independent of the choice of $B$, and moreover depends only on the $G^F$-conjugacy class of $T$ and on the orbit of $\theta$ under $(N(T)/T)^F$ where $N(T)$ is the normalizer of $T$. Note that as $\widetilde{X}_{T\subset B}$ has right $G^F\times T^F$-action, $H^i_c(\widetilde{X}_{T\subset B},\Lambda)$ is a left $G^F\times T^F$-module, and $H^i_c(X_{T\subset B},\mF_\theta)=H^i_c(\widetilde{X}_{T\subset B},\Lambda)(\theta)$, where we denote by $V(\theta)$ the maximal direct summand of $V$ on which $T^F$ acts by $\theta$.

\bde \label{geomconj}
\benu
\item Let $T,T'$ be two $F$-stable maximal tori of $G$, and $\theta,\theta'$ be characters of $T^F,T^{\prime F}$, respectively. The pairs $(T,\theta),(T',\theta')$ are said to be {\it geometrically conjugate} if the pairs $(T,\theta\circ N),(T',\theta'\circ N)$, where $N$ is the norm from $T^{F^n}$ to $T^F$ (resp.\ $T^{\prime F^n}$ to $T^{\prime F}$), are $G^{F^n}$-conjugate for some integer $n$. Here the {\it norm} $N$ for $T$ is the map $\prod_{i=0}^{n-1}F^i:T^{F^n}\ra T^F$.
\item The character $\theta$ of $T^F$ is said to be in {\it general position} if it is not fixed by any non-trivial element of $(N(T)/T)^F$.
\eenu
\ede

Here we summarize some of the main theorems in this theory:

\begin{theorem} \label{dlthms}
\benu
\item Every irreducible representation $\pi$ of $G^F$ occurs in some $R^\theta_T$, i.e.\ $\langle \pi, R_T^\theta\rangle\neq 0$ where $\langle,\rangle$ is the natural inner product on the Grothendieck group of representations of $G^F$ {\rm (\cite{DL}, Cor.\ 7.7)}.
\item If $(T,\theta)$ and $(T',\theta')$ are not geometrically conjugate, no irreducible representation of $G^F$ occurs in both $R_T^\theta$ and $R_{T'}^{\theta'}$ {\rm (\cite{DL}, Cor.\ 6.3)}. 
\item If we put, for two $F$-stable maximal tori $T,T'$,
\[N(T,T'):=\{g\in G\mid Tg=gT'\},\ \ W(T,T')^F:=T^F\setminus N(T,T')^F\]
then we have {\rm (\cite{DL}, Th.\ 6.8)}:
\[\langle R_T^\theta,R_{T'}^{\theta'} \rangle=|\{w\in W(T,T')^F\mid \theta w=w\theta'\}|.\]
Therefore, if $(T,\theta)$ is in general position, one of $\pm R_T^\theta$ (namely $(-1)^{\sigma(G)-\sigma(T)}R_T^\theta$ by {\rm \cite{DL}, Prop.\ 7.3}) is an irreducible representation of $G^F$. If moreover $T$ is not contained in any $F$-stable proper parabolic subgroup of $G$, then $(-1)^{\sigma(G)-\sigma(T)}R_T^\theta$ is a cuspidal representation.
\item The $\ol{k}$-variety $X_{T\subset B}$ is affine as long as $q$ is larger than the Coxeter number of $G$. In this case, if $\theta$ is in general position then we have $H^i_c(X_{T\subset B},\mF_\theta)=0$ for $i\neq l(w)$, where $l(w)$ is the length of the element $w\in W$ of the Weyl group $W$ such that $(B,F(B))$ is in the {\rm relative position} $w$ (i.e.\ it is in the $G$-orbit of $(B,\widetilde{w}B\widetilde{w}^{-1})$ for a representative $\widetilde{w}\in N(T)$ of $w$).
\eenu
\end{theorem}

Now we specialize to the case $G=GL_n$, and a torus $T$ associated to the element $w$ of the Weyl group corresponding to the cyclic permutation $(1,\ldots,n)$ in the symmetric group of $n$ letters, i.e.\ a torus $T$ such that $T^F\cong k_n^\times$. 

\bde
We denote by $DL$ the Deligne-Lusztig variety $\widetilde{X}_{T\subset B}$ for $G=GL_n, T^F=k_n^\times$ and $w=(1,\ldots,n)$. If $V$ is a finite dimensional representation of $GL_n(k)\times k_n^\times$ and $\theta$ is a character of $k_n^\times$, we denote by $V_\cusp$ (resp.\ $V(\theta)$) the maximal direct summand of $V$ on which $GL_n(k)$ acts by cuspidal representations (resp.\ on which $T^F$ acts by $\theta$, often considered simply as a $GL_n(k)$-representation). A character $\theta$ of $k_n^\times$ is in {\em general position} if and only if it does not factor through the norm map $k_n^\times\ra k_m^\times$ for any $m\mid n$ with $m\neq n$. We denote by $C$ the set of all characters of $k_n^\times$ in general position.
\ede

\bpr \label{DLmain}
\benu
\item {\rm (\cite{DL}, Prop.\ 7.3, Th.\ 8.3, and Cor.\ 9.9)} 
If $\theta\in C$, then:
\[ H^i_c(DL,\Lambda)(\theta) \cong
\begin{cases}
\pi_\theta\otimes \theta  & (i=n-1), \\
0 & (i\neq n-1).
\end{cases} \]
Here $\pi_\theta$ is an irreducible cuspidal representation of $GL_n(k)$ which is characterized by $\pi_\theta \otimes {\rm St}=\ind_{T(k)}^{GL_n(k)}\theta$, where ${\rm St}$ is the Steinberg representation of $GL_n(k)$. All cuspidal representations of $GL_n(k)$ arise in this way.
\item If $\theta\notin C$, then no cuspidal representation of $GL_n(k)$ occur in the cohomology groups $H^i_c(DL,\Lambda)(\theta)$ for any $i$.
\eenu
\epr

Here, (ii) follows from the slightly stronger version of the disjointness theorem (above Theorem \ref{dlthms}(ii)), which is stated only in terms of the alternating sums in \cite{DL}. It can be shown that the disjointness holds between each cohomology groups, by replacing the alternating sums in the proof of Th.\ 6.2 of \cite{DL} by each of the cohomology groups.

\bcor \label{dlcusp}
As $GL_n(k)\times k_n^\times$-representations, we have:
\[ H^i_c(DL,\Lambda)_{\rm cusp} \cong
\begin{cases}
\bigoplus_{\theta\in C}\pi_\theta\otimes \theta  & (i=n-1), \\
0 & (i\neq n-1).
\end{cases} \]
\ecor

\subsection{Computing the vanishing cycles}

Now we return to our original space $X=\Spec A$ and look at the cohomology of $X_{\ol{\eta}}$, where $\eta=\Spec \hat{K}^\ur$ is the generic point of $\Spec W$. We analyze the cohomology groups $H^i(X_{\ol{\eta}},\Lambda)$, which are finite dimensional representations of $GL_n(k)\times I_K$. 

\subsubsection{Using the first blow-up}

We start with the computation of cohomology groups using the model $Z_1$ of $X$ introduced in Section 3. By Proposition \ref{specseq}, we have 
\[ H^i(X_{\ol{\eta}},\Lambda) \cong \bH^i(Y_\Ps,R\psi\Lambda)\]
for all $i$, as $GL_n(k)\times I_K$-representations. Here we show that the cuspidal representations of $GL_n(k)$ only live in the cohomology of nearby cycle sheaves on the open subscheme $Y_\Ps^{(n)}$ of $Y_\Ps$. (Recall from the end of Section 3.4 the stratification $Y_\Ps^{(h)}\subset Y_\Ps\cong \Ps$ for $1\leq h\leq n$.)

For a $k$-rational linear subspace $N\subset \Ps$ of dimension $h-1$ and $g\in GL_n(k)$, the right $GL_n(k)$-action on $Y_\Ps$ induces $N^0\cong (Ng)^0$ and $g^*R\psi\Lambda|_{(Ng)^0}\cong R\psi\Lambda|_{N^0}$. Therefore by a standard argument (see \cite{Bo1} Lemme 13.2, \cite{DL} Prop.\ 8.2 or \cite{HT}, pp.115--117), we have the following:

\bpr \label{induced}
For a $k$-rational linear subspace $N\subset \Ps$ of dimension $h-1$, with $1\leq h\leq n$, let $P_N$ be the stabilizer of $N$, a parabolic subgroup of $GL_n(k)$. Then we have for every $i,j$:
\[ H^i_c(Y_\Ps^{(h)},R^j\psi\Lambda)\cong \ind_{P_N}^{GL_n(k)}H^i_c(N^0,R^j\psi\Lambda). \]
\epr

\bcor \label{cuspidal}
We have $H^i(X_{\ol{\eta}},\Lambda)_{\rm cusp} \cong \bH^i_c(Y_\Ps^{(n)},R\psi\Lambda)_{\rm cusp}$ for all $i$.
\ecor

\bprf
The unipotent radical $U_N$ of $P_N$ acts trivially on $N^0$, hence also on $H_c^i(N^0, R^j\psi\Lambda)$ by Proposition \ref{gensstvan}. Therefore Proposition \ref{induced} shows that $H^i_c(Y_\Ps^{(h)}, R^j\psi\Lambda)$ is parabolically induced from a representation of a Levi subgroup of $P_N$, i.e.\ $H^i_c(Y_\Ps^{(h)}, R^j\psi\Lambda)_{\rm cusp}=0$ for all $i,j$ and $h\neq n$. By using the long exact sequences:
\[\xymatrix{\cdots \ar[r] & H^i_c(Y_\Ps^{(h)}, R^j\psi\Lambda) \ar[r] & H^i(Y_\Ps^{[h]}, R^j\psi\Lambda) \ar[r] & H^i(Y_\Ps^{[h-1]}, R^j\psi\Lambda|_{Y_\Ps^{h+1}}) \ar[r] & \cdots }\]
recursively for $1\leq h\leq n-1$, we have $H^i(Y_\Ps^{[n-1]},R^j\psi\Lambda)_{\rm cusp}=0$ for any $j$, hence
\[ H^i(X_{\ol{\eta}},\Lambda)_\cusp \cong \bH^i(Y_\Ps,R\psi\Lambda)_\cusp \cong \bH^i_c(Y_\Ps^{(n)},R\psi\Lambda)_\cusp.\]
\eprf

\subsubsection{Using the generalized semistable model}

Now we make use of the generalized semistable model $Z_\st$ constructed in Section 4. By Proposition \ref{specseq}, we have 
\[ H^i(X_{\ol{\eta}},\Lambda) \cong \bH^i(Y_{\Ps,\st},R\psi\Lambda) \]
for all $i$. We will appeal to Saito's results in Section 6.1 through Proposition \ref{shgsst}. 

\bde
Let $\mG$ be the Grothendieck group of finite dimensional vector spaces over $\Lambda=\ol{\Q}_\ell$ with left action of $GL_n(k)\times I_K$. We regard the characters of $k_n^\times$ as characters of $I_K$, by composing the canonical surjection $I_K\ra k_n^\times$. 
\ede

\bde
We denote the alternating sum of the cohomology groups, regarded as elements of $\mG$, as follows:
\[H^*(X_{\ol{\eta}}):=\sum_i(-1)^i[H^i(X_{\ol{\eta}},\Lambda)],\ \ \ H^*(Y_\Ps^{(n)},R^j\psi\Lambda):=\sum_i(-1)^i[H^i_c(Y_\Ps^{(n)},R^j\psi\Lambda)].\]
\ede

\bpr \label{vancycstmodel}
\benu
\item $H^*(X_{\ol{\eta}})=H^*(Y_\Ps^{(n)},R^0\psi\Lambda)$.
\item For each degree $i$ and $\chi\in C$, we have:
\begin{align*}
H^i(X_{\ol{\eta}},\Lambda)_\cusp &\cong H^i_c(Y_\Ps^{(n)},R^0\psi\Lambda)_\cusp,\\
H^i(X_{\ol{\eta}},\Lambda)(\chi) &\cong H^i_c(Y_\Ps^{(n)},R^0\psi\Lambda)(\chi).
\end{align*}
\eenu
\epr

\bprf
Note that $R^j\psi\Lambda|_{Y_\Ps^{(n)}}=0$ for $j>0$ by Corollary \ref{corsaito}(i), hence for all $i$:
\begin{equation} \label{deg0only}
\bH^i_c(Y_\Ps^{(n)},R\psi\Lambda) \cong H^i_c(Y_\Ps^{(n)},R^0\psi\Lambda).
\end{equation}

(i) By Proposition \ref{shgsst}, we can calculate $\bH^i(Y_{\Ps,\st},R\psi\Lambda)$ on $Y_{\ol{s},\st}\subset \Sh_\st$, where we apply Corollary \ref{corsaito}(ii) to see that $\sum_i(-1)^j[R^j\psi\Lambda|_{Y_{\ol{s},\st}\setminus Y_\Ps^{(n)}}]=0$, thus:
\begin{align*}
H^*(X_{\ol{\eta}}) &= \sum_i(-1)^i[\bH^i(Y_{\Ps,\st},R\psi\Lambda)] \\
&= \sum_i(-1)^i[\bH^i_c(Y_\Ps^{(n)},R\psi\Lambda)] = H^*(Y_\Ps^{(n)},R^0\psi\Lambda).
\end{align*}

(ii) The first equality follows from Corollary \ref{cuspidal} and (\ref{deg0only}). The second is proven in an exactly similar way as in Corollary \ref{cuspidal}, except that we use Corollary \ref{corsaito}(i), instead of Proposition \ref{induced}, to see that $H^i_c(Y_{\ol{s},\st}\setminus Y_\Ps^{(n)},R^j\psi\Lambda)(\chi)=0$ for all $i,j$ (here note that multiplicities of the components intersecting $Y_\Ps$ are of the form $q^m-1$ with $m<n$, hence outside $Y_\Ps^{(n)}$ the $d$ in Corollary \ref{corsaito}(i) is strictly less than $q^n-1$ and inertia cannot act by $\chi\in C$). Use (\ref{deg0only}) to conclude the proof.
\eprf

\subsubsection{Using the model $Z_n$}

Now we have seen that the part of $H^i(X_{\ol{\eta}},\Lambda)$ we are interested in comes from the cohomology $\bH^i_c(Y_\Ps^{(n)},R\psi\Lambda)$ of nearby cycle sheaves on $Y_\Ps^{(n)}$ (it does not matter whether we regard $Y_\Ps^{(n)}$ as a subvariety of $Z_1$ or of $Z_\st$, in view of Proposition \ref{finvancyc}). Because $R^j\psi\Lambda|_{Y_\Ps^{(n)}}=0$ for $j>0$ by Corollary \ref{corsaito}(i) (or by Proposition \ref{vancycstmodel}), we can concentrate on $H^i_c(Y_\Ps^{(n)},R^0\psi\Lambda)$.

We compute this using the normalization $U_n$ and the finite etale covering $f:U_n\ra Y_\Ps^{(n)}$. By Proposition \ref{finvancyc} we see that $R^i\psi\Lambda|_{Y_\Ps^{(n)}}\cong f_*R^i\psi\Lambda|_{U_n}$ for each $i$, and as $U_n$ is the special fiber of the formally smooth $S_n$-scheme $\Spec C_n$, we have
\[R^i\psi\Lambda|_{U_n} \cong
\begin{cases}
\Lambda & (i=0)\\
0 & (i>0)
\end{cases},\ \ \ \ 
R^i\psi\Lambda|_{Y_\Ps^{(n)}} \cong
\begin{cases}
f_*\Lambda & (i=0)\\
0 & (i>0)
\end{cases}.
\]
Therefore we have a canonical $GL_n(k)\times I_K$-equivariant isomorphism
\begin{equation} \label{r0un}
H^i_c(Y_\Ps^{(n)},R^0\psi\Lambda)\cong H^i_c(U_n,\Lambda).
\end{equation}

Now we connect our result in Chapter 5 to the Deligne-Lusztig theory by comparing their explicit equations, which proves the third part of Theorem \ref{main2}:

\bpr \label{undl}
There is a $GL_n(k)\times k_n^\times$-equivariant isomorphism $U_n\cong DL$ of varieties over $\ol{k}$, where $DL$ is the Deligne-Lusztig variety defined in Section 6.2.
\epr

\bprf
This is readily seen by comparing our equation in Proposition \ref{equ2} and the $GL_n(k)\times k_n^\times$-action in Proposition \ref{actionun} with the equation and group actions of $DL$ given in Section 2.2 of \cite{DL}.
\eprf

Therefore we can invoke the Deligne-Lusztig theory to derive our main theorem. Let us denote the alternating sum of cohomology groups of $DL$ as:
\[H^*_c(DL):=\sum_i(-1)^i[H^i_c(DL,\Lambda)]=\sum_\theta R_T^\theta\otimes[\theta]\in \mG.\]
For a character $\chi\in C$ of $T^F\cong k_n^\times$ in general position, considered as a character of $I_K$, recall the corresponding cuspidal representation $\pi_\chi$, defined in Proposition \ref{DLmain}.

\bth \label{cohomology}
\benu
\item $H^*(X_{\ol{\eta}})=H^*_c(DL)$.
\item For each degree $i$ and $\chi\in C$, we have:
\begin{align*}
H^i(X_{\ol{\eta}},\Lambda)_{\rm cusp} &\cong
\begin{cases}
\bigoplus_{\chi\in C}\pi_\chi\otimes \chi & (i=n-1), \\
0 & (i\neq n-1),
\end{cases} \\
H^i(X_{\ol{\eta}},\Lambda)(\chi) &\cong
\begin{cases}
\pi_\chi\otimes \chi & (i=n-1), \\
0 & (i\neq n-1).
\end{cases}
\end{align*}
\eenu
\ethm

\bprf
We combine the isomorphism (\ref{r0un}) and Proposition \ref{undl} with the Proposition \ref{vancycstmodel}, and use the Deligne-Lusztig theory (Proposition \ref{DLmain} and Corollary \ref{dlcusp}) for (ii).
\eprf

{\bf Note added in proof}: In the proof of Proposition 6.10, 
we need a little more argument to prove that $U_N$ acts trivially 
on $R^j\psi\Lambda$. We use $Z_\st$ to apply
Proposition 6.3. Although the first isomorphism of Proposition 6.3 
holds only etale locally, it shows that the canonical morphism
$\Lambda \rightarrow R^0\psi\Lambda|_{Y_J^0}$ is an isomorphism if $d=1$.
For general $d$, by Proposition 6.2, the sheaf $R^0\psi\Lambda$ is a push 
forward from a $d=1$ situation, namely the normalization of the base change 
to tamely ramified extension of $W$ of degree $d$ (similar to what
is done in \S5.1), which restricts to a finite etale covering of degree $d$ on $Y_J^0$.
Thus $U_N$, being a $p$-group, acts trivially on $R^0\psi\Lambda|_{Y_J^0}$. 
In the second isomorphism of Proposition 6.3, the group $U_N$ can 
only act on the index set $J$, but $J$ is a partial flag of linear 
subspaces of $\Ps$ containing $N$, and $U_N$ fixes each element of $J$.

\end{document}